\documentclass[11pt]{article}
\usepackage{amsmath,amstext,amssymb,amsfonts,amsthm}
\usepackage{mathrsfs}
\usepackage{subfigure}
\usepackage{epsfig}
\usepackage[ruled]{algorithm2e}

\usepackage{graphicx}
\usepackage{accents}
\usepackage{color}
\usepackage{bm}
\usepackage[export]{adjustbox}
\usepackage{epstopdf}

\newtheorem{lem}{Lemma}[section]
\newtheorem{thm}{Theorem}[section]

\newtheorem{Def}{Definiton}[section]

\numberwithin{equation}{section}
\numberwithin{figure}{section}

\newcommand{\bn}{\boldsymbol{n}}
\newcommand{\bw}{\boldsymbol{w}}
\newcommand{\bt}{\boldsymbol{t}}
\newcommand{\bx}{\boldsymbol{x}}
\newcommand{\bbe}{\boldsymbol{e}}
\newcommand{\bu}{\boldsymbol{u}}
\newcommand{\bz}{\boldsymbol{z}}
\newcommand{\by}{\boldsymbol{y}}
\newcommand{\bg}{\boldsymbol{g}}
\newcommand{\ba}{\boldsymbol{a}}
\newcommand{\bv}{\boldsymbol{v}}

\newcommand{\bV}{\boldsymbol{V}}
\newcommand{\bF}{\boldsymbol{F}}
\newcommand{\bE}{{\boldsymbol{E}}}
\newcommand{\bJ}{{\boldsymbol{J}}}
\newcommand{\bX}{{\boldsymbol{X}}}
\newcommand{\bL}{{\boldsymbol{L}}}
\newcommand{\bW}{{\boldsymbol{W}}}
\newcommand{\bQ}{{\boldsymbol{Q}}}
\newcommand{\bH}{{\boldsymbol{H}}}
\newcommand{\bP}{{\boldsymbol{\Phi}}}

\newcommand{\om}{\omega}
\newcommand{\pa}{\partial}
\newcommand{\la}{\langle}
\newcommand{\ra}{\rangle}

\newcommand{\al}{\alpha}

\newcommand{\ep}{\varepsilon}
\newcommand{\na}{\nabla}
\newcommand{\vp}{\varphi}
\newcommand{\ga}{\gamma}
\newcommand{\Ga}{\Gamma}
\newcommand{\Om}{\Omega}
\newcommand{\Si}{\Sigma}
\newcommand{\de}{\delta}

\newcommand{\De}{\Delta}

\newcommand{\lam}{\lambda}

\newcommand{\vep}{\varepsilon}
\newcommand{\lj}{[{\hskip -1.5pt}[}
\newcommand{\rj}{]{\hskip -1.5pt} ]}
\newcommand{\ljj}{\Big[{\hskip -3pt}\Big[}
\newcommand{\rjj}{\Big]{\hskip -3pt}\Big]}

\newcommand{\dist}{\mathrm{dist}}

\newcommand{\cT}{{\cal T}}
\newcommand{\cM}{{\cal M}}

\newcommand{\cF}{{\cal F}}

\newcommand{\cO}{{\cal O}}

\renewcommand{\i}{\mathbf{i}}
\newcommand{\R}{{\mathbb{R}}}

\newcommand{\C}{{\mathbb{C}}}
\newcommand{\A}{{\mathbb{A}}}

\renewcommand{\div}{\mathrm{div}\,}
\newcommand{\curl}{\mathrm{curl}}

\newcommand{\be}{\begin{eqnarray}}
\newcommand{\ee}{\end{eqnarray}}
\newcommand{\beq}{\begin{equation}}
\newcommand{\eeq}{\end{equation}}
\newcommand{\ben}{\begin{eqnarray*}}
\newcommand{\een}{\end{eqnarray*}}
\newcommand{\nn}{\nonumber}

\newcommand{\ddhat}[1]{\widehat{\vphantom{\rule{3pt}{8.5pt}}\smash{\widehat{#1}}}}

\newcommand{\rev}[1]{{\color{black}{#1}}}
\newcommand{\revto}[1]{{\color{black}{#1}}}
\newcommand{\revt}[1]{{\color{black}{#1}}}

\title{A high order unfitted finite element method for time-harmonic Maxwell interface problems
\thanks{This work is supported in part by China National Key Technologies R\&D Program under the grant 2019YFA0709602 and China Natural Science Foundation under the grant 118311061, 12288201, 12101609.}}
\author{Zhiming Chen\thanks{LSEC, Institute of Computational Mathematics,
Academy of Mathematics and System Sciences and School of Mathematical Science, University of
Chinese Academy of Sciences, Chinese Academy of Sciences,
Beijing 100190, China. ({\tt zmchen@lsec.cc.ac.cn}).}
\and Ke Li \thanks{Department of Basic, Information Engineering University, Zhengzhou 450001, China. ({\tt like@lsec.cc.ac.cn}).}
\and Maohui Lyu\thanks{School of Sciences, Beijing University of Posts and Telecommunications, Beijing 100876, China. ({\tt mlyu@bupt.edu.cn}).}
\and Xueshuang Xiang\thanks{Qian Xuesen Laboratory of Space Technology,
    China Academy of Space Technology, Beijing 100194, China. ({\tt xiangxs@lsec.cc.ac.cn}).}}

\begin{document}

\maketitle

\begin{center}
\small
\begin{minipage}{0.9\textwidth}
\textbf{Abstract.}
We propose a high order unfitted finite element method for solving time-harmonic Maxwell interface problems. The unfitted finite element method is based on a mixed formulation in the discontinuous Galerkin framework on a Cartesian mesh with possible hanging nodes. The $H^2$ regularity of the solution to Maxwell interface problems with $C^2$ interfaces in each subdomain is proved. Practical \rev{interface-resolving} mesh conditions are introduced under which the $hp$ inverse estimates on three-dimensional curved domains are proved.  Stability and $hp$ a priori error estimate of the  unfitted finite element method are proved. Numerical results are included to illustrate the performance of the method.

\medskip

\end{minipage}
\end{center}

\setlength{\parindent}{2em}
%%%%%%%%%%%%%%%%%%%%%%%%%%%%%%%%%%%%%%%%%%%%%%%%%%%
\section{Introduction}
%%%%%%%%%%%%%%%%%%%%%%%%%%%%%%%%%%%%%%%%%%%%%%%%%%%

Let $\Om\subset \R^3$ be a bounded domain with a Lipschitz boundary $\Sigma$. We consider in this paper the following time-harmonic Maxwell interface problem
\be
& &\na\times(\mu^{-1}\na\times \bE)-k^2\vep\bE=\bJ,\ \ \div(\vep\bE)=0\ \ \ \ \mbox{in }\Om,\label{p1}\\
%& &\ \ \mbox{in }\Om,\label{p2}\\
& &\lj\bE\times\bn\rj_\Ga=0,\ \ \lj(\mu^{-1}\na\times\bE)\times\bn\rj_\Ga=0,\ \ \lj\vep\bE\cdot\bn\rj_\Ga=0\ \ \ \ \mbox{on }\Ga,\label{p3}\\
& &\bE\times\bn=\bg\times\bn\ \ \ \ \mbox{on }\Sigma,\label{p4}
\ee
where $\bJ\in\bL^2(\Om)$ with $\div\bJ=0$ in $\Om$ and $\bg\times\bn\in\bH^{\rev{3}/2}(\Sigma)$. Here and 
throughout the paper, for any Banach space $X$, we denote $\bX= X^3$ and $\|\cdot\|_X$ both the norms of $X$ and $\bX$. 

We assume the domain $\Om$ is divided by a $C^2$ interface $\Ga$ into two subdomains so that $\Om=\Om_1\cup\Ga\cup\Om_2$ and 
\rev{$\Om_1$ is strictly included in} $\Om$.
For simplicity, we assume the relative permeability $\mu$ and the relative permittivity $\vep$ are piecewise \rev{constants}
$\mu=\mu_1\chi_{\Om_1}+\mu_2\chi_{\Om_2}$, $\vep=\vep_1\chi_{\Om_1}+\vep_2\chi_{\Om_2}$,
where for $i=1,2$, $\chi_{\Om_i}$ is the characteristic function of $\Om_i$, and $\mu_i,\vep_i$ are positive constants. $k=\omega\sqrt{\vep_0\mu_0}$ is the wave number of the vacuum with $\om>0$ the angular frequency and $\mu_0,\vep_0$ the permeability and permittivity of the vacuum. With this notation, $\bJ=\i k\mu_0\bJ_a$ with $\bJ_a$ being the applied current density. We denote \rev{by} $\bn$ both the unit outer normal to $\Om_1$ on $\Ga$ and the unit outer normal to $\Om$ on $\Sigma$. $\lj v\rj_\Ga:=v|_{\Om_1}-v|_{\Om_2}$ stands for the jump of a function $v$ across the interface $\Ga$. 

The existence and uniqueness of the weak solution to the problem \eqref{p1}-\eqref{p4} are well studied in the literature (see, e.g., \cite{Hi02}). The singularity and regularity of the solution with smooth $\mu,\vep$ in polyhedral and smooth domains are considered in \cite{Co99, Co10}.
The singularity of the solution of the Maxwell interface problems with polyhedral interfaces is studied in \cite{Co00}. To the best of the authors' knowledge, the $H^2$ interface regularity of the solution to the Maxwell interface problem with smooth interfaces is missing in the literature. In this paper, we first prove the $H^2$ regularity of the solution to \eqref{p1}-\eqref{p4} in each subdomain $\Om_1,\Om_2$. Our proof is based on the $\bH(\curl)$-coercive Maxwell equations which is different from the $\bH(\curl)\cap\bH(\div)$-coercive Maxwell equations used in \cite{Co99, Co00, Co10}. This new regularity result will be used in our finite element convergence analysis based on the Schatz argument in dealing with the indefiniteness of the time-harmonic Maxwell equations. 

There exists a large literature on finite element methods for solving the time-harmonic Maxwell equations after the seminar work
\cite{Ne80}. We refer to \cite{Hi02, Mo03} for the study of the $\bH(\curl)$-conforming $h$-methods, \cite{De08, Me19} for the $hp$-methods, and \cite{Pe03, Ho04, Ho05} for the discontinuous Galerkin (DG) methods. The common assumptions in these studies are that the domains are polyhedral and the material interfaces are piecewise flat so that conforming tetrahedral or hexahedral meshes can be used. Much less studies have been devoted to finite element methods solving Maxwell equations on domains with curved boundary. We refer to \cite{Du90, Hi12, Me18} for body-fitted finite element methods, \cite{Bu10} for the isogeometric analysis, and \cite{Ch09} for the low order unfitted finite element method. We remark that the design of body-fitted high-order finite element methods depends on nonlinear element transforms from the reference element to the elements with curved boundary \cite{Be89, Me18}. It may be challenging to satisfy the conditions on the nonlinear element transforms which depend on the geometry of the interface or boundary in practical applications.

The original motivation of unfitted finite element methods in the DG framework is to release the time-consuming work of constructing shape regular meshes for domains with complex geometry. It turns out that the unfitted finite element methods also provide a natural way to design high-order methods without resorting to the nonlinear element transforms. Since the seminal work \cite{Ha02} for elliptic interface problems, considerable progress of the unfitted finite element methods has been made in the literature \cite{Burman10, Wu, Massjung, Jo13, Badia18, Wang, Liu20, Ch20}. The small cut cell problem, that is, the intersection of the domain and the elements can be arbitrarily small or anisotropic, can be solved by appropriate techniques of stabilization \cite{Burman10, Wu, Massjung} or merging small cut cells with surrounding large elements \cite{Jo13, Badia18, BurmanHHO21, Ch20}. We refer to \cite{Ch20, Ch22} for further references and other approaches of unfitted finite element methods.

In \cite{Ch20} an adaptive high-order unfitted finite element method in two dimension is proposed for elliptic interface problems in which the $hp$ a priori and a posteriori error estimates are derived based on novel $hp$ domain inverse estimates and the concept of interface deviation. The interface deviation is a measure that quantifies the resolution of the geometry by the mesh. In \cite{Ch22}, for any $C^2$ interface, a reliable algorithm is constructed to merge small interface elements with their surrounding elements to generate an induced mesh whose elements are large with respect to both domains, which solves the small cut cell problem.  It is also shown in 
\cite{Ch22} that the exponential growth of the condition number of the stiffness matrix on the finite element approximation order, which is observed in the literature (e.g., \cite{de17, Ch22}), can be controlled by reducing the interface deviation. Therefore, arbitrarily high order unfitted finite element methods on automatically generated Cartesian meshes for solving elliptic interface problems can be achieved for arbitrarily shaped $C^2$ smooth interfaces.

The main purpose of this paper is to extend the high order unfitted finite element method for two-dimensional elliptic interface problems in \cite{Ch20} to solve the time-harmonic Maxwell equations \eqref{p1}-\eqref{p4}. We characterize and quantify the mesh resolution of the geometry in two steps. We first introduce the concept of proper intersection of the interface and boundary to the elements and the concept of large element in three dimension, which allow us to show that each large element is a union of strongly shape regular polyhedrons in the sense that the polyhedron is a union of shape regular tetrahedrons in the classical sense. Then we introduce the concept of \rev{the} interface deviation in three dimension and prove the $hp$ inverse trace estimate on curved domains which plays an important role in studying the unfitted finite element method in this paper. 

One of the key difficulties in the finite element convergence analysis for the indefinite time-harmonic Maxwell equations is the discretization of the
divergence free condition $\div(\vep\bE)=0$. In the $\bH(\curl)$-conforming finite element method, this condition is implemented implicitly and the finite element convergence analysis depends on the fact that any discrete divergence free finite element function can be approximated by a continuous divergence free function (see, e.g., \cite{Hi02, Mo03}). This property, which depends crucially on the construction of conforming edges elements, unfortunately, 
is not available for unfitted finite element functions which can not be conforming across the curved interface. 
To overcome the difficulty, we propose to use the mixed formulation to explicitly enforce the discrete divergence free 
condition and penalize
the jump of the normal component of the discrete electric field in the $H^{1/2}$ norm, which implies to penalize $\div(\vep\bE)=0$
in the $L^2$ norm, see Lemma \ref{lem:4.3} below. The stability and $hp$-error estimates are proved in the asymptotic range, that is, sufficiently small mesh sizes. 

The layout of the paper is as follows. In section 2 we introduce some notation and prove the $H^2$ regularity of the problem \eqref{p1}-\eqref{p4} in each subdomain $\Om_1$, $\Om_2$. In section 3 we introduce the unfitted finite element method for the time-harmonic
Maxwell equations. In section 4 we prove the stability and $hp$-a priori error estimate for our unfitted finite element method. In section 5 we report some numerical results. 

%%%%%%%%%%%%%%%%%%%%%%%%%%%%%%%%%%%%%%%%%%%%%%%%%%%
\section{The Maxwell interface problem} 
%%%%%%%%%%%%%%%%%%%%%%%%%%%%%%%%%%%%%%%%%%%%%%%%%%%
We first recall some notation. For any Lipschitz domain $D$ in $\R^3$ with boundary $\pa D$ whose unit normal is denoted by $\bn$, the space
$\bH(\curl;D)=\{\bv\in\bL^2(D):\na\times\bv\in\bL^2(D)\}$
is a Hilbert space under the graph norm. For any $\bv\in\bH(\curl;D)$, its tangential trace is defined as $\ga_\tau(\bv)=\bv\times\bn$ on $\pa D$. It is known that (see Buffa et al. \cite{Bu02}) $\ga_\tau:\bH(\curl;D)\to\bH^{-1/2}(\div_{\pa D};\pa D)$ is bounded and surjective, where
\ben
\bH^{-1/2}(\div_{\pa D};\pa D)=\{\boldsymbol{\lam}\in \bV_\pi(\pa D)':\div_{\pa D}(\boldsymbol{\lam})\in\bH^{-1/2}(\pa D)\}.
\een
Here $\bV_\pi(\pa D)'$ is the dual space of $\bV_\pi(\pa D)=\{\bv_T=\bn\times\bv\times\bn:\bv\in\bH^{1/2}(\pa D)\}$ and $\div_{\pa D}:\bV_\ga(\pa D)\to H^{-1/2}(\pa D)$ is the dual operator of the surface gradient $\na_T:H^{1/2}(\pa D)\to \bV_\ga(\pa D)'$, where $\bV_\ga(\pa D)=\{
\bv\times\bn:\bv\in\bH^{1/2}(\pa D)\}$. It is known that 
$\div_{\pa D}(\bv\times\bn)=(\na\times\bv)\cdot\bn\in\bH^{-1/2}(\pa D)$ for any $\bv\in\bH(\curl;D)$. We denote $\bH_0(\curl;D):=\{\bv\in\bH(\curl;D):\bv\times\bn=0\ \ \mbox{on }\pa D\}$.  

The following two lemmas will be used in the proof of the $H^2$ regularity up to the interface of the solution to time-harmonic Maxwell equations.

\begin{lem}\label{lem:2.2}
Let $D$ be a bounded Lipschitz domain, $\bF\in\bH^1(D)'\revto{\cap\bH_0(\curl;D)'}$, and $\div\bF\in H^{-1}(D)$. Then there exists a constant $C$ depending only on the domain $D$ such that
$\|\bF\|_{\bH_0(\curl;D)'}\le C(\|\bF\|_{H^1(D)'}+\|\div\bF\|_{H^{-1}(D)})$.
\end{lem}

\begin{proof}
By the Birman-Solomyak regular decomposition theorem \cite{Bi87}, Hiptmair \cite[Lemma 2.4]{Hi02}, any $\bv\in\bH_0(\curl;D)$ can be splitted as $\bv=\bv_s+\na\psi$ for
some $\bv_s\in\bH^1(D)$, $\psi\in H^1_0(D)$ such that $\|\bv_s\|_{H^1(D)}+\|\psi\|_{H^1(D)}\le C\|\bv\|_{\bH(\curl;D)}$. 
\revto{Since $\psi\in H^1_0(D)$, for any $\de>0$, there exists $\psi_\de\in C^\infty_0(D)$ such that $\psi_\de\to\psi$ in $H^1_0(D)$ as $\de\to 0$, which also implies $\na\psi_\de\to\na\psi$ in $\bH_0(\curl;D)$ as $\de\to 0$. By the definition of the derivatives of the distribution, $(\div\bF)(\psi_\de)=-\bF(\na\psi_\de)$. Now since $\bF\in\bH_0(\curl;D)'$ and $\div\bF\in H^{-1}(D)$, we obtain $\bF(\na\psi)=-(\div\bF)(\psi)$ by letting $\de\to 0$.}
Thus
\ben
\|\bF\|_{\bH_0(\curl;D)'}&=&\sup_{\bv\in \bH_0(\curl;D)}\frac{|\revto{\bF(\bv)}|}{\|\bv\|_{\bH(\curl;D)}}\nn\\
&\le&C\sup_{\bv_s\in\bH^1(D),\psi\in H^1_0(D)}\frac{|\revto{\bF(\bv_s+\na\psi)}|}{\|\bv_s\|_{H^1(D)}+\|\psi\|_{H^1(D)}}\nn\\
&\le&C(\|\bF\|_{H^1(D)'}+\|\div\bF\|_{H^{-1}(D)}).
\een
This completes the proof.
\end{proof}

\begin{lem}\label{lem:2.3}
Let $D$ be a bounded Lipschitz domain which is divided by a $C^{1,1}$ surface $\Ga_D$ into two parts $D_1,D_2$.
Let $\mathbb{B}\in \R^{3\times 3}$ be a symmetric matrix with elements in $C^{0,1}(\bar D_1\cup\bar D_2)$ such that $a_-|\bx|^2\le\mathbb{B}\bx\cdot\bx\le a_+|\bx|^2\ \ \forall\bx\in\R^3$ a.e. in $D$ for some constants $a_-,a_+>0$.
Let $f\in L^2(D), g\in H^{1/2}(\Ga_D)$, and $\bv\in\bH(\curl;D)$ satisfy  
\beq
\div(\mathbb{B}\bv)=f\ \ \mbox{in }D_1\cup D_2,\ \ \lj\mathbb{B}\bv\cdot\bn\rj_{\Ga_D}=g\ \ \mbox{on }\Ga_D. \label{xx1}
\eeq
Then for any subdomain \revt{$\mathcal{O}$ strictly included in $D$}, $\mathcal{O}_1=\mathcal{O}\cap D_1$, $\mathcal{O}_2=\mathcal{O}\cap D_2$, we have $\bv\in\bH^1(\mathcal{O}_1\cup \mathcal{O}_2)$, and 
$\|\bv\|_{H^1(\mathcal{O}_1\cup \mathcal{O}_2)}\le C(\|\bv\|_{\bH(\curl;D)}+\|f\|_{L^2(D)}+\|g\|_{H^{1/2}(\Ga_D)})$,
where the constant $C$ depends on the domains $\mathcal{O}$, $D$, and the coefficient matrix $\mathbb{B}$.
\end{lem} 

\begin{proof} The argument is standard (see, e.g., Hiptmair et al. \cite[Lemma 3.2]{Hi11}). We sketch a proof for the sake of completeness. By the regular decomposition theorem, $\bv=\bv_s+\na\psi$ for some $\bv_s\in\bH^1(D)$, $\psi\in H^1(D)$ such that
$\|\bv_s\|_{H^1(D)}+\|\psi\|_{H^1(D)}\le C\|\bv\|_{\bH(\curl;D)}$.
It is clear from \eqref{xx1} that $\psi\in H^1(D)$ satisfies
\ben
& &\div(\mathbb{B}\na\psi)=f-\div(\mathbb{B}\bv_s)\ \ \mbox{in }D_1\cup D_2,\\
& &\lj\psi\rj_{\Ga_D}=0,\ \ \ljj\mathbb{B}\frac{\pa\psi}{\pa\bn}\rjj_{\Ga_D}=g-\lj\mathbb{B}\bv_s\cdot\bn\rj_{\Ga_D}\ \ \ \ \mbox{on }\Ga_D.
\een
By the regularity theorem of elliptic interface problems (see, e.g., McLean \cite[Theorem 4.20]{Mc00}), $\psi\in H^2(\cO_1\cup \cO_2)$ and 
\ben
\|\psi\|_{H^2(\cO_1\cup \cO_2)}
&\le&C(\|\psi\|_{H^1(D)}+\|f-\div(\mathbb{B}\bv_s)\|_{L^2(D)}+\|g-\lj\mathbb{B}\bv_s\cdot\bn\rj_{\Ga_D}\|_{H^{1/2}(\Ga_D)})\\
&\le&C(\|\bv\|_{\bH(\curl;D)}+\|f\|_{L^2(D)}+\|g\|_{H^{1/2}(\Ga_D)}).
\een
This completes the proof.
\end{proof}

The weak formulation of the problem \rev{\eqref{p1}-\eqref{p4}} is to find ${\bE}\in\bH(\curl;\Om)$ such that ${\bE}\times\bn=\bg\times\bn$ on $\Sigma$, and
\beq\label{aa3}
(\mu^{-1}\na\times\bE,\na\times\bv)-k^2(\varepsilon{\bE},\bv)=(\bJ,\bv)\ \ \forall\bv\in\bH_0(\curl;\Om),
\eeq
where $(\cdot,\cdot)$ denotes the inner product on $\bL^2(\Om)$. Notice that \eqref{aa3} implies $\div(\vep\bE)=0$ in $\Om$
by taking $\bv=\na\psi$ for any $\psi\in H^1_0(\Om)$. It is known that except a denumerable wave numbers $0<k_1<k_2<\cdots$ of nonzero Maxwell eigenvalues tending to $\infty$, the problem \eqref{aa3} exists a unique solution (see, e.g., \cite[\S4.1-\S4.2]{Hi02}). In this paper, we will assume $k$ is not equal to one of these eigenvalues and thus the problem \eqref{aa3} has a unique solution $\bE\in\bH(\curl;\Om)$ which satisfies \cite[(96)]{Hi02}
\beq\label{aa7}
\|{\bE}\|_{\bH(\curl;\Om)}\le C(\|\bJ\|_{L^2(\Om)}+\|\bg\times\bn\|_{\bH^{-1/2}(\div_\Sigma;\Sigma)}).
\eeq

The following $H^1$ regularity of the solution of \eqref{aa3} can be proved by the same argument used in the proof of Lemma \ref{lem:2.3}. Here we omit the details.

\begin{lem}\label{lem:2.5}
Let the interface $\Ga$ and the boundary $\Sigma$ be $C^{1,1}$, $\bJ\in\bL^2(\Om),\bg\in\bH^{1/2}(\Sigma)$ satisfy $\div\bJ=0$ in $\Om$. Then the solution of the problem \eqref{aa3} satisfies $\bE\in\bH^1(\Om_1\cup\Om_2)$ and 
$\|\bE\|_{H^1(\Om_1\cup\Om_2)}\le C(\|\bJ\|_{L^2(\Om)}+\|\bg\|_{H^{1/2}(\Sigma)})$,
where the constant $C$ depends on the domain $\Om$ and the coefficients $k,\mu,\vep$.
\end{lem}

The following theorem is the main result of this section. It will be used in our finite element convergence analysis in section 4. 

\begin{thm}\label{thm:2.1}
Let the interface $\Ga$ and the boundary $\Sigma$ be $C^2$, $\bJ\in\bL^2(\Om)$ with $\div\bJ=0$ in $\Om$, and $\bg\in\bH^{3/2}(\Sigma)$. Then the solution of \eqref{aa3} has the regularity $\bE\in\bH^2(\Om_1\cup\Om_2)$ and satisfies
$\|\bE\|_{H^2(\Om_1\cup\Om_2)}
\le C_{\rm reg}(\|\bJ\|_{L^2(\Om)}+\|\bg\|_{H^{3/2}(\Sigma)})$,
where the constant $C_{\rm reg}>0$ depends on the domain $\Om$, the interface $\Ga$, the boundary $\Sigma$, and the coefficients $k,\vep,\mu$.
\end{thm}

\begin{proof} \rev{The proof follows the idea to prove the $H^2$ regularity for elliptic equations (see, e.g., Evans \cite[Chapter 6]{Evans98}, \cite[Chapter 4]{Mc00}) based on selecting smooth cut-off functions to localize the equation and flatting the interface and the boundary by means of a suitable change of variables.} 
Let $\{B(\bx_i,r_i)\}^n_{i=1}$, where $\bx_i\in\Ga$, $i=1,\cdots,n$, be a union of balls inside $\Om$ such that $\cO_\Ga=\cup^n_{i=1}B(\bx_i,r_i/2)$ covers $\Ga$, and $\{B(\by_i,d_i)\}^m_{i=1}$, where $\by_i\in\Sigma$, $i=1,\cdots,m$, be a union of balls not intersecting with the interface $\Ga$ such that $\cO_\Sigma=\cup^m_{i=1}B(\by_i,d_i/2)$ covers $\Sigma$. Let $\cO_i=\Om_i\backslash(\bar{\cO}_\Ga\cup \bar{\cO}_\Sigma)$, $i=1,2$. Clearly, the union of
$\cO_\Ga,\cO_\Sigma$, and $\cO_1\cup\cO_2$ covers $\Om$. The proof is divided into three parts.

$1^\circ$ Interior regularity. Since $\rev{\div(\vep\bE)}=0$ in $\Om_1\cup\Om_2$, by Lemma \ref{lem:2.5}, $\bE\in\bH^1(\Om_1\cup\Om_2)$ and it satisfies
\ben
-\Delta\bE=\mu\bJ+k^2\mu\vep\bE\ \ \ \ \mbox{in }\Om_1\cup\Om_2.
\een
By the interior regularity of elliptic equations (see, e.g., \cite[Theorem 4.16]{Mc00}) and Lemma \ref{lem:2.5}, we know that $\bE\in\bH^2(\cO_i)$, $i=1,2$, and
\be
\hskip-0.6cm\|\bE\|_{H^2(\cO_i)}&\le&C(\|\bE\|_{H^1(\Om_i)}+\|\bJ+k^2\mu\vep\bE\|_{L^2(\Om_i)})\le C(\|\bJ\|_{L^2(\Om)}+\|\bg\|_{H^{1/2}(\Sigma)}).\label{w1}
\ee

$2^\circ$ Interface regularity. For $i=1,\cdots,n$, let $B=B(\bx_i,r_i)$ and $\hat B$ be the unit ball in $\R^3$. Let $\bP:\hat B\to B$ be the $C^2$, one-to-one mapping such that $B_1=\Om_1\cap B=\bP(\hat B_1), B_2=\Om_2\cap B=\bP(\hat B_2)$, and $\Ga\cap B=\bP(\hat \Ga)$, where $\hat B_1=\{\hat \bx\in\hat B:\hat x_3<0\}$, $\hat B_2=\{\hat \bx\in\hat B:\hat x_3>0\}$, and $\hat\Ga=\{\hat \bx\in\hat B:\hat x_3=0\}$. By the Piola transform (see, e.g., Monk \cite[\S 3.9]{Mo03}), we have
\beq
\na=D\bP^{-T}\hat\na,\ \ \na\cdot=J^{-1}\hat\na\cdot(JD\bP^{-1}),\ \ \na\times=J^{-1}D\bP\hat\na\times D\bP^T,\label{a5}
\eeq
where $D\bP$ is the gradient matrix and $J={\rm det}(D\bP)$. Moreover, the unit normal $\bn$ to $\Ga$ and the surface area $ds$ of $\Ga$ satisfy (see, e.g., Hofmann et al. \cite{Ho07})
\beq\label{a4}
\bn\circ\bP={D\bP^{-T}\hat\bn}/{|D\bP^{-T}\hat\bn|},\ \ ds=|J|\,|D\bP^{-T}\hat\bn|\,d\hat s,
\eeq
where $\hat\bn=(0,0,1)^T$ and $d\hat s$ is the surface area of $\hat\Ga$.

Let $\chi\in C^\infty_0(B)$ be the cut-off function such that $0\le\chi\le 1$ in $B$ and $\chi=1$ in $B(\bx_i,r_i/2)$. Denote $\bu=\chi\bE$. Then we have
\be
& &\na\times(\mu^{-1}\na\times\bu)=\bF,\ \ \div(\vep\bu)=\na\chi\cdot(\vep\bE)\ \ \mbox{in }\Om_1\cup\Om_2,\label{a6a}\\
& &\lj\bu\times\bn\rj_\Ga=0,\ \ \lj(\mu^{-1}\na\times\bu)\times\bn\rj_\Ga=\bg_1\times\bn,\ \ \lj\vep\bu\cdot\bn\rj_\Ga=0\ \ \mbox{on }\Ga,\label{a6b}
\ee
where $\bg_1\times\bn=\lj(\mu^{-1}\na\chi\times\bE)\times\bn\rj_{\Ga}$ on $\Ga$, and
\beq\label{a7}
\bF=\bF_1+\na\times(\mu^{-1}\na\chi\times\bE),\ \ \bF_1=\chi\bJ+k^2\vep\chi\bE+\na\chi\times(\mu^{-1}\na\times\bE)\ \ \mbox{in }\Om.
\eeq 
Now for any function $\bv:B\to\C^3$, we denote $\hat\bv=D\bP^T(\bv\circ\bP)$. By \eqref{a5}-\eqref{a4} and the following vector identity \cite{Ho07}
\beq\label{v1}
{\mathbb D}\boldsymbol{a}\times {\mathbb D}\boldsymbol{b}=\det(\mathbb{D}){\mathbb D}^{-T}\boldsymbol{a}\times\boldsymbol{b}\ \ \ \ \forall {\mathbb D}\in\C^{3\times 3},\ \ \forall\boldsymbol{a, b}\in\C^3,
\eeq
we deduce from \eqref{a6a}-\eqref{a6b} that 
\be
& &\hat\na\times(\hat\mu^{-1}\mathbb{A}\hat\na\times\hat\bu)=\mathbb{A}^{-1}\hat\bF,\ \  \hat\na\cdot(\hat\vep\mathbb{A}^{-1}\hat\bu)=\mathbb{A}^{-1}\hat\na\hat\chi\cdot(\hat\vep\hat\bE)\ \ \mbox{in }\hat B_1\cup\hat B_2,\label{a8a}\\
& &\lj\hat\bu\times\hat\bn\rj_{\hat\Ga}=0,\ \ \lj(\hat\mu^{-1}\mathbb{A}\hat\na\times\hat\bu)\times\hat\bn\rj_{\hat\Ga}=\hat\bg_1\times\hat\bn,\ \ \lj\hat\vep\mathbb{A}^{-1}\hat\bu\cdot\hat\bn\rj_{\hat\Ga}=0\ \ \mbox{on }\hat\Ga,\label{a8b}
\ee
where $\hat\mu=\,\mu\circ\bP$, $\hat\vep=\vep\circ\bP$, $\hat\chi=\chi\circ\Phi$, $\mathbb{A}=J^{-1}D\bP^TD\bP$, and 
\beq\label{a9}
\mathbb{A}^{-1}\hat\bF=\mathbb{A}^{-1}\hat\bF_1+\hat\na\times(\hat\mu^{-1}\mathbb{A}(\hat\na\hat\chi\times\hat\bE)),\ \ \hat\bg_1\times\hat\bn=\lj(\hat\mu^{-1}\mathbb{A}(\hat\na\hat\chi\times\hat\bE))\times\hat\bn\rj_{\hat\Ga}.
\eeq 

Since $\hat\bn=(0,0,1)^T$, $\lj\hat\bu\times\hat\bn\rj_{\hat\Ga}=0$ implies that $\lj\hat u_1\rj_{\hat\Ga}=\lj\hat u_2\rj_{\hat\Ga}=0$ on $\hat\Ga$. For any $\de>0$ sufficiently small, we consider the difference quotients $\De_l^\de$, $l=1,2$, defined as
\ben
\De^\de_l\hat\bu(\hat{\bx})=(\hat\bu(\hat\bx+\de\bbe_l)-\hat\bu(\hat\bx))/\de.
\een
By \eqref{a8a}-\eqref{a8b} we have
\ben
& &\hat\na\cdot(\hat\vep\mathbb{A}^{-1}_\de\Delta^\de_l\hat\bu)=-\hat\na\cdot(\hat\vep(\Delta^\de_l\mathbb{A}^{-1})\hat\bu)+\Delta^\de_l(\mathbb{A}^{-1}\hat\na\hat\chi\cdot(\hat\vep\hat\bE))\ \ \mbox{in }\hat B_1\cup\hat B_2,\\
& &\lj\hat\vep\mathbb{A}^{-1}_\de\Delta^\de_l\hat\bu\cdot\hat\bn\rj_{\hat\Ga}=-\lj\hat\vep(\Delta_l^\de\mathbb{A}^{-1})\hat\bu\cdot\hat\bn\rj_{\hat\Ga}\ \ \mbox{on }\hat\Ga,
\een
where $\mathbb{A}_\de=\mathbb{A}(\hat\bx+\de\bbe_l)$. Since $\bP^{-1}(B(\bx_i,r_i/2)\cap\Om_k)\subset\hat B_k$, $k=1,2$, by Lemma \ref{lem:2.3}, 
\be
& &\|\De^\de_l\hat\bu\|_{H^1(\bP^{-1}(B(\bx_i,r_i/2)\cap(\Om_1\cup\Om_2)))}\nn\\
&\le&C(\|\De^\de_l\hat\bu\|_{\bH(\curl;\hat B)}+\|\hat\bu\|_{H^1(\hat B_1\cup \hat B_2)}+\|\hat\bE\|_{H^1(\hat B_1\cup\hat B_2)})\nn\\
&\le&C(\|\hat\na\times\Delta^\de_l\hat\bu\|_{L^2(\hat B)}+\|\bE\|_{H^1(B_1\cup B_2)}).\label{a11}
\ee
Next by \eqref{a8a}-\eqref{a8b} we have
\be
& &\hat\na\times(\hat\mu^{-1}\A_\de\hat\na\times(\De^\de_l\hat\bu))=\hat\bF'\ \ \mbox{in }\hat B_1\cup\hat B_2,\label{a10a}\\
& &\lj(\De^\de_l\hat\bu)\times\hat\bn\rj_{\hat\Ga}=0,\ \ \lj(\hat\mu^{-1}\mathbb{A}_\de\hat\na\times(\De^\de_l\hat\bu))\times\hat\bn\rj_{\hat\Ga}=\hat\bg_1'\times\hat\bn\ \ \mbox{on }\hat\Ga,\label{a10b}
\ee
where $\hat\bg'_1\times\hat\bn=(\De^\de_l\hat\bg_1)\times\hat\bn-\lj(\hat\mu^{-1}(\De^\de_l\mathbb{A})\hat\na\times\hat\bu)\times\hat\bn\rj_{\hat\Ga}$, and
\beq\label{z5}
\hat\bF'=\De^\de_l(\mathbb{A}^{-1}\hat\bF)-\hat\na\times(\hat\mu^{-1}(\De^\de_l\mathbb{A})\hat\na\times\hat\bu).
\eeq 
Since $(\De^\de_l\hat\bu)\times\hat\bn=0$ on $\pa\hat B$, by testing \eqref{a10a} with $\De^\de_l\hat\bu$ we obatin
\beq
\|\hat\na\times\De^\de_l\hat\bu\|_{L^2(\hat B)}\le C\,\sup_{\hat\bv\in\bH_0(\curl;\hat B)}\frac{|(\hat\bF',\hat\bv)_{\hat B}+\la\hat\bg_1'\times\hat\bn,\hat\bv_T\ra_{\hat\Ga}|}{\|\hat\bv\|_{\bH_0(\curl;\hat B)}}.\label{a12}
\eeq
From the definition of $\hat\bF'$ and $\hat\bg_1'\times\hat\bn$, we obtain by doing integration by parts that
\ben
(\hat\bF',\hat\bv)_{\hat B}+\la\hat\bg_1'\times\hat\bn,\hat\bv_T\ra_{\hat\Ga}
&=&(\De^\de_l(\mathbb{A}^{-1}\hat\bF),\hat\bv)_{\hat B}+\la\De^\de_l\hat\bg_1\times\hat\bn,\hat\bv_T\ra_{\hat\Ga}\nn\\
& &\ -\,(\hat\mu^{-1}(\De^\de_l\A)(\hat\na\times\hat\bu),\hat\na\times\hat\bv)_{\hat B}.
\een
Since $\hat\bF$ is compactly supported in $\hat B$ and $\hat\bg_1\times\hat\bn$ is compactly supported on $\hat B\cap\hat\Ga$, by the definition of $\mathbb{A}^{-1}\hat\bF$ and $\hat\bg_1\times\hat\bn$ in \eqref{a9} and integration by parts, we have then, for sufficiently small $\de$,
\ben
& &(\De^\de_l(\mathbb{A}^{-1}\hat\bF),\hat\bv)_{\hat B}+\la\De^\de_l\hat\bg_1\times\hat\bn,\hat\bv_T\ra_{\hat\Ga}\\
&=&-(\mathbb{A}^{-1}\hat\bF,\De^{-\de}_l\hat\bv)_{\hat B}-\la\hat\bg_1\times\hat\bn,\De^{-\de}_l\hat\bv_T\ra_{\hat\Ga}\\
&=&-(\mathbb{A}^{-1}\hat\bF_1,\De^{-\de}_l\hat\bv)_{\hat B}-(\hat\mu^{-1}\mathbb{A}(\hat\na\hat\chi\times\hat\bE),\hat\na\times\De^{-\de}_l\hat\bv)_{\hat B}\\
&=&-(\De^\de_l(\mathbb{A}^{-1}\hat\bF_1),\hat\bv)_{\hat B}
-(\hat\mu^{-1}\De^\de_l(\A(\hat\na\hat\chi\times\hat\bE)),\hat\na\times\hat\bv)_{\hat B}.
\een
By Lemma \ref{lem:2.2} 
\be
\|\De^\de_l(\mathbb{A}^{-1}\hat\bF_1)\|_{\bH_0(\curl;\hat B)'}&\le&C(\|\De^\de_l(\mathbb{A}^{-1}\hat\bF_1)\|_{(H^{1}(\hat B))'}+\|\hat\na\cdot(\De^\de_l(\mathbb{A}^{-1}\hat\bF_1))\|_{H^{-1}(\hat B)})\nn\\
&\le&C\|\mathbb{A}^{-1}\hat\bF_1\|_{H(\div;\hat B)}.\label{r1}
\ee
Thus, by \eqref{a7} and $\div\bF_1=0$ in $\Om$, \eqref{a12} can be bounded as
\be\label{c4}
\|\hat\na\times\De^\de_l\hat\bu\|_{L^2(\hat B)}&\le&C(\|\mathbb{A}^{-1}\hat\bF_1\|_{H(\div;\hat B)}+\|\hat\bE\|_{H^1(\hat B_1\cup\hat B_2)})\nn\\
&\le&C(\|\bF_1\|_{H(\div;B)}+\|\bE\|_{H^1(B_1\cup B_2)})\nn\\
&\le&C(\|\bJ\|_{L^2(B)}+\|\bE\|_{H^1(B_1\cup B_2)}).
\ee
Combining \eqref{a11} with \eqref{c4}, we obtain
\ben
\|\De^\de_l\hat\bu\|_{H^1(\bP^{-1}(B(\bx_i,r_i/2)\cap (\Om_1\cup\Om_2)))}\le C(\|\bJ\|_{L^2(B)}+\|\bE\|_{H^1(B_1\cup B_2)}).
\een
This implies, by letting $\de\to 0$, for $l=1,2$,
\beq
\|\pa_{\hat x_l}\hat\bu\|_{H^1(\bP^{-1}(B(\bx_i,r_i/2)\cap (\Om_1\cup\Om_2)))}
\le C(\|\bJ\|_{L^2(B)}+\|\bE\|_{H^1(B_1\cup B_2)}).\label{c5}
\eeq
Finally, it follows from \eqref{a8a} that
\ben
-\hat\na\cdot(\hat\vep \A^{-1}\hat\bu)=-\mathbb{A}^{-1}\hat\na\hat\chi\cdot(\hat\vep\hat\bE)\ \ \ \ \mbox{in }\hat B_1\cup\hat B_2.
\een
Notice that $\mathbb{A}^{-1}=J D\bP^{-1} D\bP^{-T}$ whose elements $(\mathbb{A}^{-1})_{kl}=J\ba^T_k\ba_l$, $k,l=1,2,3$, where $\ba_1,\ba_2,\ba_3\in\R^3$ are column vectors of $D\bP^{-T}$. Obviously, $(\mathbb{A}^{-1})_{33}=J|\ba_3|^2\ge a_0$ for some constant $a_0>0$. Thus by differentiating the equation in $\hat x_3$ we obtain 
\ben
\left\|\frac{\pa^2\hat u_3}{\pa \hat x_3^2}\right\|_{L^2(\bP^{-1}(B(\bx_i,r_i/2)\cap (\Om_1\cup\Om_2)))}&\le&
C\sum^2_{l=1}\|\pa_{\hat x_l}\hat\bu\|_{H^1(\bP^{-1}(B(\bx_i,r_i/2)\cap (\Om_1\cup\Om_2)))}\\
&&\,+\,C(\|\hat \bE\|_{H^1(\hat B_1\cup\hat B _2)}+\|\hat\bu\|_{H^1(\hat B _1\cup\hat B _2)}).
\een
Therefore, it follows from \eqref{c5} that
\ben
\|\hat\bu\|_{H^2(\Phi^{-1}(B(\bx_i,r_i/2)\cap (\Om_1\cup\Om_2))}\le C(\|\bJ\|_{L^2(B)}+\|\bE\|_{H^1(B_1\cup B_2)}),
\een
and thus $\|\bu\|_{H^2(B(\bx_i,r_i/2)\cap (\Om_1\cup\Om_2)))}\le C(\|\bJ\|_{L^2(B)}+\|\bE\|_{H^1(B_1\cup B_2)})$ for any $i=1,\cdots,n$. This yields by Lemma \ref{lem:2.5} that
\beq
\|\bE\|_{H^2(\mathcal{O}_\Ga\cap (\Om_1\cup\Om_2))}\le C(\|\bJ\|_{L^2(\Om)}+\|\bg\|_{H^{1/2}(\Sigma)}).\label{w2}
\eeq

$3^\circ$ Boundary regularity.  Let $\bu_{\bg}\in\bH^2(\Om)$ be the lifting of $\bg\in\bH^{3/2}(\Sigma)$ such that $\|\bu_{\bg}\|_{H^2(\Om)}\le C\|\bg\|_{H^{3/2}(\Sigma)}$. Then $\bu=\bE-\bu_{\bg}$ satisfies $\bu\times\bn=0$ on $\Sigma$, and
\ben
\na\times(\mu^{-1}\na\times\bu)-k^2\vep\bu=\bJ'\ \ \ \ \mbox{in }\Om,
\een
where $\bJ'=\bJ-\na\times(\mu^{-1}\na\times\bu_{\bg})+k^2\vep\bu_{\bg}$ in $\Om$. Similar to the argument in the step $2^\circ$, for $i=1,\cdots,m$, we consider $\bu_i=\chi_i\bu$, where $\chi_i\in C^\infty_0(B(\by_i,d_i))$ is the cut-off function such that $0\le\chi_i\le 1$, $\chi_i=1$ in $B(\by_i,d_i/2)$, and use the $H^2$ regularity of the solution to 
time harmonic Maxwell equations in $C^2$ domains in Costabel et al. \cite[Theorems 4.5.3 and 3.4.5]{Co10} to obtain
\beq\label{w3}
\|\bE\|_{H^2(\cO_\Sigma\cap\Om)}\le C (\|\bJ\|_{L^2(\Om)}+\|\bg\|_{H^{3/2}(\Sigma)}).
\eeq
The theorem follows now from \eqref{w1}, \eqref{w2}, and \eqref{w3}. 
\end{proof}

\revto{We remark that it is crucial to use Lemma \ref{lem:2.2} in \eqref{r1} to bound $\De^\de_l(\mathbb{A}^{-1}\hat{\bF}_1)$ in the $\bH_0(\curl;\hat B)'$ norm. If one simply bounds $\De^\de_l(\mathbb{A}^{-1}\hat{\bF}_1)$ in the $L^2$ norm, one will have to require $\bJ\in\bH^1(B)$ in \eqref{c4}. The regularity bound in Theorem \ref{thm:2.1} will become $\|\bE\|_{H^2(\Om_1\cup\Om_2)}\le C_{\rm reg}(\|\bJ\|_{H^1(\Om)}+\|\bg\|_{H^{3/2}(\Om)})$, which is not sufficient for us to use the Schatz argument in the proof of Theorem \ref{thm:4.2}.}

To conclude this section, we introduce the mixed formulation for \eqref{p1}-\eqref{p4} to be used in this paper. Let $\vp\in H^1_0(\Om)$ be the Lagrangian multiplier for the constraint $\div(\vep\bE)=0$ in $\Om$ in \eqref{p1}. The mixed formulation is to find $(\bE,\vp)\in\bH(\curl;\Om)\times H^1_0(\Om)$ such that $\bE\times\bn=\bg\times\bn$ on $\Sigma$, and
\be
& &(\mu^{-1}\na\times\bE,\na\times\bv)-k^2(\vep\bE,\bv)-(\vep\na\vp,\bv)
=(\bJ,\bv)\ \ \forall \bv\in \bH_0(\curl;\Om),\label{a1}\\
& &(\vep\bE,\na\psi)=0\ \ \ \ \forall\psi\in H^1_0(\Om).\label{a2}
\ee
It is easy to see that if the wave number $k$ is not equal to the Maxwell eigenvalues and $\bE$ is the unique solution of \eqref{aa3}, then the mixed problem \eqref{a1}-\eqref{a2} has a unique solution $(\bE,\vp)$ with $\vp=0$ in $\Om$.

%%%%%%%%%%%%%%%%%%%%%%%%%%%%%%%%%%%%%%%%%%%%%%%%%%%
\section{The unfitted finite element method}
%%%%%%%%%%%%%%%%%%%%%%%%%%%%%%%%%%%%%%%%%%%%%%%%%%%

In this section we first introduce some notation on the finite element mesh and prove the crucial inverse trace estimate on curved domains by extending 
the idea in Chen et al. \cite{Ch20} for the two-dimensional unfitted finite element method. Then we introduce the unfitted finite element method for the mixed formulation of the time-harmonic Maxwell interface problem \eqref{a1}-\eqref{a2}.

\subsection{Notation and the inverse trace estimate}

Let $\cT$ be a mesh consisting of right hexahedrons with possible hanging nodes that covers $\Om$. For each $K\in\cT^\Ga=\{K\in\cT:K\cap\Ga\not=\emptyset\}$ and $K\in\cT^\Si=\{K\in\cT:K\cap\Si\not=\emptyset\}$, we assume the interface $\Ga$ or the boundary $\Sigma$ intersects $K$ properly in the following sense.

\begin{Def}\label{def:3.0} {\rm (Proper intersection)}
\rev{The intersection of the interface $\Ga$ or the boundary $\Sigma$ with an element $K$ is called the proper intersection if $\Ga$ or $\Sigma$ intersects each (open) edge of $K$ at most once, and if $\Ga$ or $\Sigma$ intersects a face $F$ of $K$, then $\Ga$ or $\Si$ intersects the edges of $F$ at most twice at different (open) edges.}
\end{Def}

\rev{If $\Ga$ or $\Si$ intersects an element at some vertex, we regard $\Ga$ or $\Si$ as intersecting one of the three edges originated from $A$ at some point very close to $A$. If $\Ga$ or $\Si$ is tangent to an edge or a face of an element, we regard $\Ga$ or $\Si$ as being very close to the edge or the face but not intersecting with the edge or the face.}

\begin{figure}[h]
	\centering
	\includegraphics[width=1\textwidth]{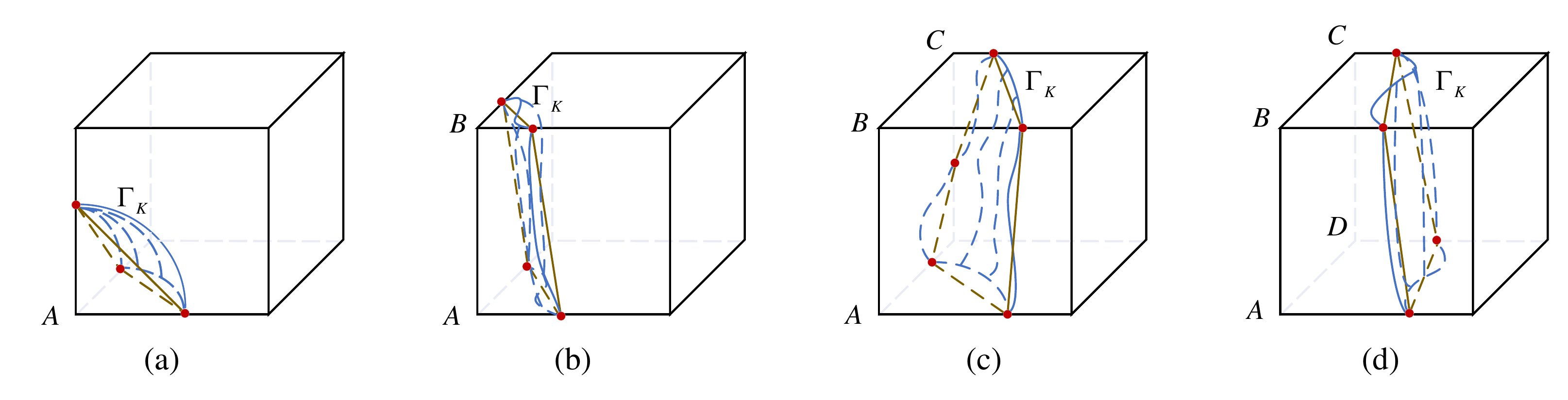}
	\caption{Illustration of proper intersections of the interface $\Ga$ and an element $K$. From (a) to (d), $K$ has $1,2,3,4$ vertices in one of the domains $\Om_i$, $i=1,2$.}
	\label{fig:3.1}
\end{figure}

\begin{figure}[h]
	\centering%\vskip-1.5cm
	\includegraphics[width=0.9\textwidth]{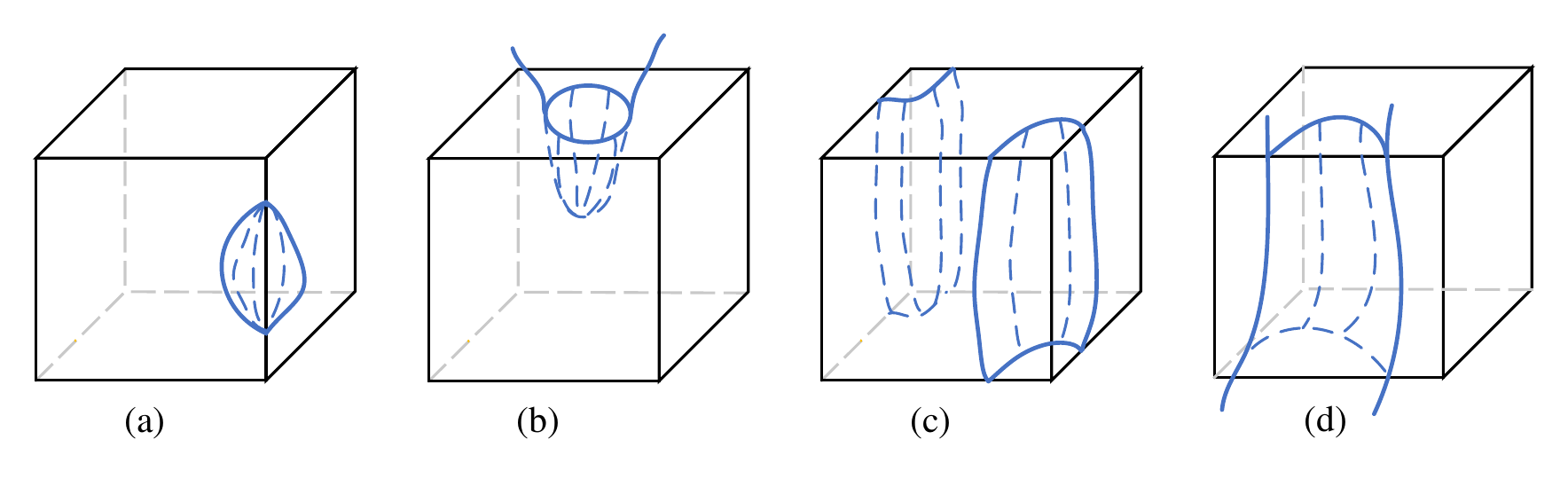}%\vskip-2cm
	\caption{Examples of improper intersections of the interface $\Ga$ and an element $K$, for which local refinements 
	are required to resolve the interface.}
	\label{fig:3.2}
\end{figure}

\revto{It is clear that when $\Ga$ and $\Si$ have proper intersections with all elements in $\cT^\Ga\cup\cT^\Si$, if two vertices $A,B$ of an element are inside $\Om_i$, then the edge $AB$ is included in $\Om_i$, $i=1,2$.}
Fig.\ref{fig:3.1} shows some of the situations that are allowed by Definition \ref{def:3.0}. In the situations shown in Fig.\ref{fig:3.2}, local refinements are needed to resolve the interface or the boundary. 

\begin{Def}\label{def:3.1} {\rm (Large element)}
	For $i=1,2$, an element $K\in\cT$ is called a large element with respect to $\Om_i$ if $K\subset\Om_i$;
	or $K\in\cT^\Ga\cup\cT^\Si$ for which there exists a \rev{fixed} $\de_0\in (0,1/2)$ such that $|e\cap\Om_i|\ge\de_0|e|$ for each edge $e$ of $K$ having nonempty intersection with $\Om_i$. An element $K\in\cT^\Ga\cup\cT^\Sigma$ is called a large element if $K\in\cT^\Ga$ is large with respect to both $\Om_1,\Om_2$ or $K\in\cT^\Si$ is large with respect to $\Om_2$. An element $K\in\cT^\Ga\cup\cT^\Sigma$ is called a small element if it is not a large element.
\end{Def}

Following \cite{Ch20} we make the following assumption on the finite element mesh when $K\in\cT^\Ga\cup\cT^\Sigma$ is a small element. 

\medskip
{\bf Assumption (H1)} 
For each $K\in\cT^\Ga$ or $K\in\cT^\Si$, there exists a \rev{cuboid} macro-element $N(K)$ which is a union of $K$ and its surrounding element (or elements) such that (i) $\Ga$ or $\Sigma$ has proper intersection with $N(K)$; and (ii) $N(K)$ is large with respect to both $\Om_1,\Om_2$ or is large with respect to $\Om_2$. We assume $h_{N(K)}\le C_0h_K$ for some fixed constant $C_0$. Here $h_K$ is the diameter of $K$.
\medskip

\begin{figure}[h]
	\centering%\vskip-1.5cm
	\includegraphics[width=0.65\textwidth]{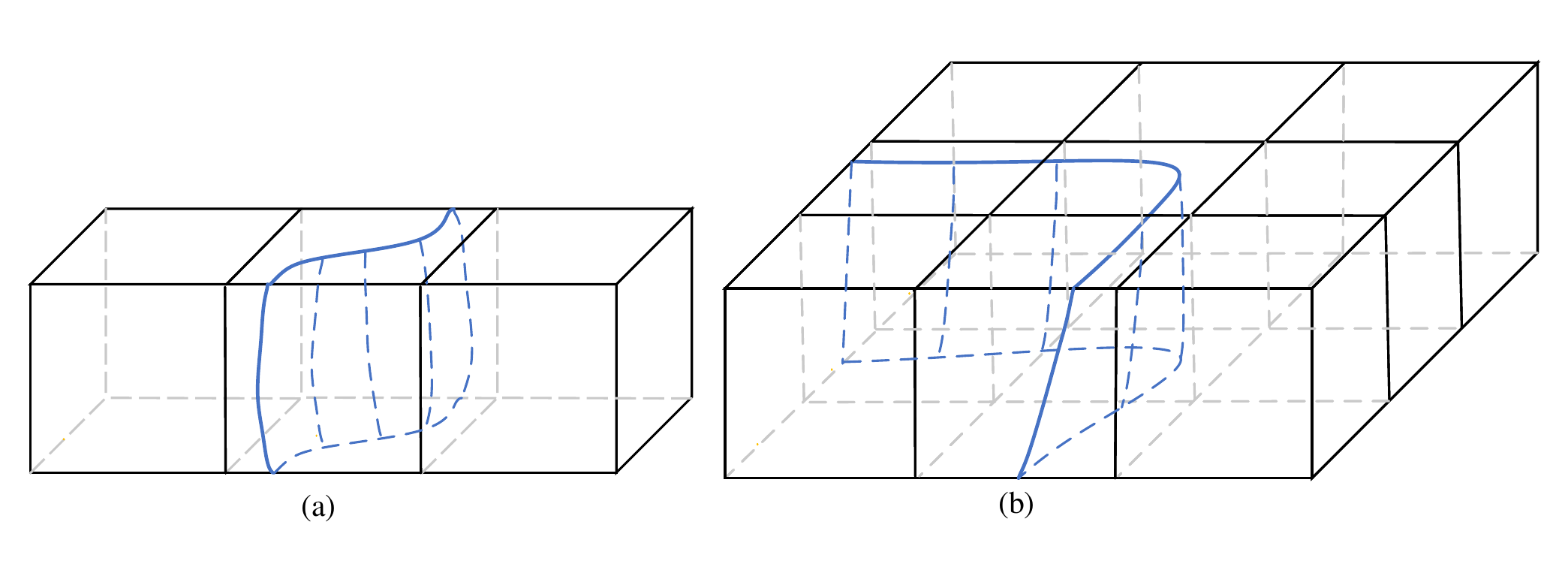}%\vskip-2cm
	\caption{Examples of merging a small element $K$ with 2 or 8 of its neighboring elements to form a large element.}
	\label{fig:3new}
\end{figure}

Inspired by Johansson and Larson \cite{Jo13} in the setting of
a fictitious boundary DG method for elliptic equations, one possible way to satisfy the assumption (H1) is to locally refine the surrounding elements $K'$ of a small element $K\in\cT^\Ga$ so that the elements $K'$ are of the same size as $K$ and $K'$ are completely included in $\Om_1$ or $\Om_2$. Then we can define $N(K)$ as the union of $K$ and those elements $K'$. A similar construction can be made for small boundary elements $K\in\cT^\Si$. \rev{Fig.\ref{fig:3new} illustrates a small element merged with 2 or 8 of its neighboring elements to form a macro-element. When a small element $K$ is of the form shown in Fig.\ref{fig:3.1}(a), the macro-element may include all 26 elements having non-empty intersection with the element $K$.}
We refer to Chen and Liu \cite{Ch22} for a reliable algorithm to merge small interface elements with their surrounding elements to automatically generate a finite element mesh whose elements are large with respect to both domains $\Om_1,\Om_2$ for any two-dimensional smooth interfaces. It is expected such a merging algorithm can also be constructed for the three dimensional smooth interfaces and will be pursued in a future work.

In the following, we will always set $N(K)=K$ if $K\in\cT^\Ga\cup\cT^\Si$ is a large element. Thus $\cM=\{N(K):K\in\cT^\Ga\cup\cT^\Si\}\cup\{K\in\cT: K\subset\Om_i, i=1,2, K\not\subset N(K') \mbox{ for some $K'\in\cT^\Ga\cup\cT^\Si$}\}$ is also a mesh that covers $\Om$ consisting of right hexahedrons. We will call $\cM$ the induced mesh of $\cT$ and write $\cM={\rm Induced }\,(\cT)$. We set $\cM^\Ga=\{K\in\cM:K\cap\Ga\not=\emptyset\}$, $\cM^\Si=\{K\in\cM:K\cap\Si\not=\emptyset\}$. 

For any $K\in\cM^\Ga$, let $\Ga_K=K\cap\Ga$ and denote \rev{by} $\mathcal{V}_K$ the set of intersection points of $\Ga$ and the edges of $K$. Similarly, for any $K\in\cM^\Si$, let $\Si_K=K\cap\Si$ and denote \rev{by} $\mathcal{V}_K'$ the set of intersection points of $\Si$ and the edges of $K$.

We call a polyhedron is {\it strongly} shape regular if it is the union of shape regualr tetrahedrons in the classical sense of Ciarlet \cite[P.132]{Ci78}. The following lemma shows that if $K\in\cT^\Ga$ is a large element, it is the union of two strongly shape regular polyhedrons. In this paper, all polyhedrons are considered to be open. We call a polyhedron $D$ is the union of polyhedrons $T_i$, $i=1,\cdots,I$, if $\bar D=\cup^I_{i=1}\bar T_i$.

\begin{figure}[h]
	\centering%\vskip-0.8cm
	\includegraphics[width=0.43\textwidth]{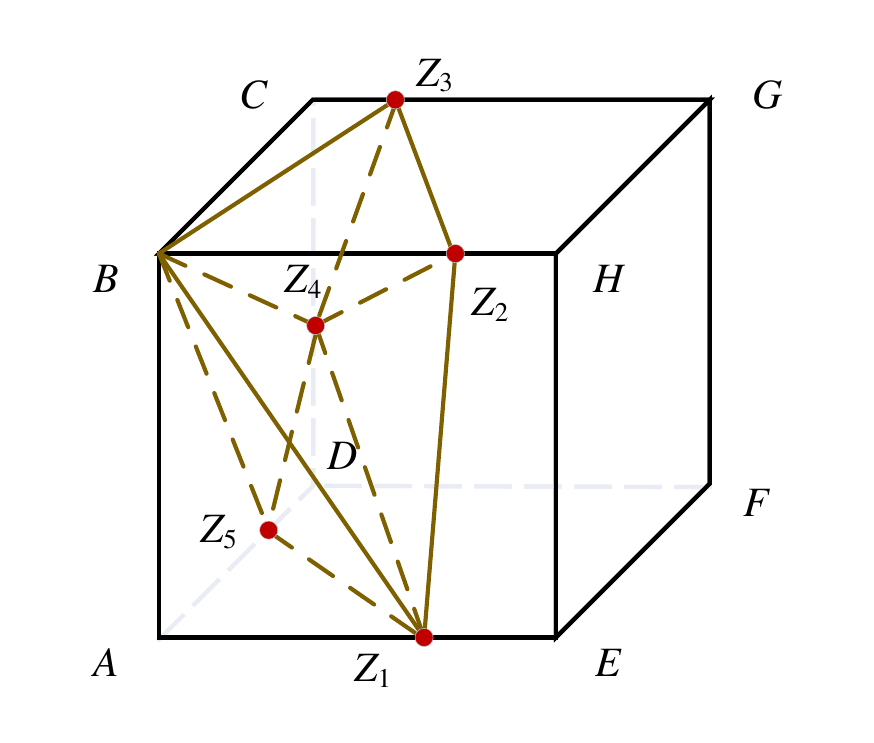}
	\includegraphics[width=0.45\textwidth]{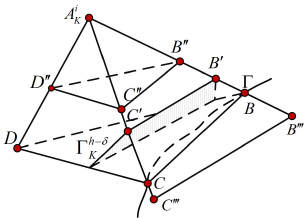}%\vskip-0.8cm
	\caption{Left: The figure used in the proof of Lemma \ref{lem:3.1new}. Right: The figure used in the proof of Lemma \ref{lem:3.4}. The tetrahedron with at most two intersection points of $\Ga$ and the edges of $K$.}
	\label{fig:3.3}
\end{figure}

\begin{lem}\label{lem:3.1new}
If $K\in\rev{\cM^\Ga}$, then $\bar K=\bar K_1^h\cup \bar K_2^h$, where for $i=1,2$, $K_i^h$ is the union of tetrahedrons $T_{ij}$, $j=1,\cdots,m_i$, $m_i\ge 1$, whose vertices are the vertices of $K$ inside $\Om_i$ and the points in $\mathcal{V}_K$. The tetrahedrons $T_{ij}$, $j=1,\cdots,m_i$, have a common vertex $A_K^i$ which is a vertex of $K$ inside $\Om_i$. Moreover, $T_{ij}$, $j=1,\cdots,m_i$, are shape regular in the sense that the radius of the inscribed ball of $T_{ij}$ is bounded below by $c_0h_K$ for some constant $c_0>0$ depending only on $\de_0$ in Definition \ref{def:3.1}.  
\end{lem}

\begin{proof} First notice that $K\in\cM^\Ga$ is a large element. We only prove the case when $K$ has three vertices $A,B,C$ inside $\Om_1$. The other cases can be proved analogously. Without loss of generality, we assume $A=(0,0,0)^T, B=(0,0,h_3)^T, C=(0,h_2,h_3)^T$, where $h_i=h_i(K)$ is the length of the side of $K$ in the $x_i$ direction, $i=1,2,3$, see Fig.\ref{fig:3.3} (left). Since if two vertices inside $\Om_1$, the edges connecting the vertices is also in $\Om_1$, there are five intersecting points $Z_j$, $j=1,\cdots,5$, of $\Ga_K$ and the edges of $K$. The coordinates of the points are $Z_1=(a_1,0,0)^T$, $Z_2=(a_2,0,h_3)^T$, $Z_3=(a_3,h_2,h_3)^T$, $Z_4=(0,h_2,c_4)^T$, $Z_5=(0,b_5,0)^T$. By Definition \ref{def:3.1}, $\de_0\le a_jh_1^{-1}\le 1-\de_0$, $j=1,2,3$, $\de_0\le c_4h_3^{-1}\le 1-\de_0$, and $\de_0\le b_5h_2^{-1}\le 1-\de_0$.

We recall that for any tetrahedron $T$, the radius $\rho_T$ of the inscribed ball of $T$ is $\rho_T=3|T|/|\pa T|$. It is easy to see that the tetrahedron $T_{11}$ with vertices $B,A,Z_1,Z_5$ and $T_{12}$ with vertices $B,Z_4,Z_3,C$ satisfy $|T_{11}|, |T_{12}|\ge \frac 16\de_0^2|K|$. The volume of the tetrahedron $T_{13}$ with vertices $B,Z_4,Z_1,Z_5$ is
\ben
|T_{13}|=\frac 16\left|\begin{array}{ccc}
0 & h_2 & c_4-h_3 \\
a_1 & 0 & -h_3 \\
0 & b_5 & -h_3 \\
\end{array}\right|=\frac 16a_1\left[h_2h_3+b_5(c_4-h_3)\right]=\frac 16a_1[(h_2-b_5)h_3+b_5c_4],
\een
which yields $|T_{13}|\ge \frac 16(\de_0+\de_0^2)\de_0|K|$. Similarly, the volume of the tetrahedron $T_{14}$ with vertices $B,Z_1,Z_2,Z_4$ is $|T_{14}|=\frac 16a_2h_2h_3\ge\frac 16\de_0|K|$, and the volume of the tetrahedron $T_{15}$ with vertices $B,Z_4,Z_2,Z_3$ is $|T_{15}|=\frac 16a_2h_2(h_3-c_4)\ge\frac 16\de_0^2|K|$. On the other hand, it is easy to see that for $j=1,\cdots,5$, $|\pa T_{1j}|\le Ch_K^2$ for some constant $C$ depending only on $\de_0$. Thus $\bar K_1^h=\cup^5_{j=1}\bar T_{1j}$ with $T_{1j}$, $j=1,\cdots,5$, being shape regular tetrahedrons. 

Let $K_2^h=K\backslash\bar K_1^h$. One can show that $K_2^h$ is a union of shape regular tetrahedrons with the common vertex $F$ in Fig.\ref{fig:3.3} (left). Here we omit the details. This completes the proof.
\end{proof}

We remark that in the proof it is important to have $B$ as the common vertex of the tetrahedrons $T_{1j}$, $j=1,\cdots,5$. If one chooses $A$ as the common vertex, the volume of the tetrahedron $T_0$ with vertices $A,Z_2,Z_3,Z_4$ is $|T_0|=\frac 16h_2[a_2(c_4-h_3)+a_3h_3]$, which may not have a uniform lower bound with respect to the positions of $Z_2,Z_3,Z_4$.

The following lemma can be proved by the same argument as that in Lemma \ref{lem:3.1new}.

\begin{lem}\label{lem:3.2new}
If $K\in\rev{\cM^\Sigma}$ is a large element, then there is a polyhedron $K^h\subset K$, $K^h$ is the union of tetrahedrons $T_j$, $j=1,\cdots,m$, $m\ge 1$, whose vertices are the vertices of $K$ inside $\Om$ and the points in $\mathcal{V}_K'$. The tetrahedrons $T_j$, $j=1,\cdots,m$, have a common vertex $A'_K$ which is a vertex of $K$ inside $\Om$. Moreover,  $T_{j}$, $j=1,\cdots,m$, are shape regular in the sense that the radius of the inscribed ball of $T_{j}$ is bounded below by $c_0h_K$ for some constant $c_0>0$ depending only on $\de_0$ in Definition \ref{def:3.1}. 
\end{lem}

We remark that the construction of the strongly shape regular polyhedrons $K_i^h$ in Lemma \ref{lem:3.1new} and $K^h$ in Lemma \ref{lem:3.2new} may not be unique. For example, in the case when $K$ has two or four vertices in $\Om_i$, each one of the vertices $K$ in $\Om_i$ can be chosen as the common vertex $A^i_K$ of the shape regular tetrahedrons. In the following, we always fix one construction of the strongly shape regular polydrons $K_1^h,K_2^h$ for $K\in\cM^\Ga$, and $K^h$ for $K\in\cM^\Sigma$.

Now for any $K\in\cM^\Ga$, we want to approximate $\Ga_K$ with a flat face $\Ga_K^h$ not far from $\Ga_K$. Since $K$ is a large element, by Lemma \ref{lem:3.1new} we know that $K$ is the union of two strongly regular polyhedrons $K_1^h,K_2^h$. It is clear that $\pa K_1^h\cap \pa K_2^h=\cup^{J_K}_{j=1}f_j$, where $\{f_j\}^{J_K}_{j=1}$ are triangles with vertices in $\mathcal{V}_K$ which may not be coplanar. Let $F_j$ be the plane that $f_j$ lies in, $j=1,\cdots,J_K$. We define $\Ga_K^h$ as one of the planes among $\{F_j\}^{J_K}_{j=1}$ such that $\dist_{\rm H}(\Ga_K,\Ga_K^h)=\min_{1\le j\le J_K}\dist_{\rm H}(\Ga_K,F_j)$, where $\dist_{\rm H}(\Ga_1,\Ga_2)=\max_{\bx\in\Ga_1}(\min_{\by\in\Ga_2}|\bx-\by|)$. Roughly speaking, ${\rm dist}_{\rm H}(\Ga_1,\Ga_2)$ measures how far $\Ga_1$ deviates from $\Ga_2$. Let $A_K^i$ be the common vertex in $\Om_i$ of the tetrahedron defined in Lemma \ref{lem:3.1new}.

Similarly, for $K\in\cM^\Sigma$, let $K^h$ be the strongly shape regular polyhedron defined in Lemma \ref{lem:3.2new}. Then $\pa K^h\backslash\pa K=\cup^{J'_K}_{j=1}f_j'$, where $f_j'$ are the triangles with vertices in $\mathcal{V}_K'$, $1\le j\le J'_K$. We define $\Si_K^h$ as one of the planes $\{F_j'\}^{J'_K}_{j=1}$ such that $\dist_{\rm H}(\Sigma_K,\Sigma^h_K)=\min_{1\le j\le J'_K}\dist_{\rm H}(\Si_K,F_j')$, where $F_j'$ is the plane that $f_j'$ lies in, $j=1,\cdots,J'_K$. Let $A_K'$ be the vertex of $K$ in $\Om$ defined in Lemma \ref{lem:3.2new}.

The following definition extends that in \cite{Ch20} for two-dimensional interfaces.

\begin{Def}\label{def:1.2} {\rm (Interface and boundary deviation)}
The interface deviation $\eta_K$ for $K\in\cM^\Ga$ and the boundary deviation $\eta_K$ for $K\in\cM^\Si$ are defined as 
	\ben
	\eta_K=\max_{i=1,2}\frac{\dist_{\rm H}(\Ga_K,\Ga_K^h)}{\dist(A_K^i,\Ga_K^h)}\ \ \forall K\in\cM^\Ga,\ \ \ \ 
	\eta_K=\frac{\dist_{\rm H}(\Si_K,\Si_K^h)}{\dist(A_K',\Si_K^h)}\ \ \forall K\in\cM^\Si,
	\een
where $\dist(A,\Ga_1)=\min_{\by\in\Ga_1}|A-\by|$ is the distance of the point $A$ to the set $\Ga_1\subset\R^3$.
\end{Def}

\begin{lem}\label{lem:3.0}
If the interface $\Ga$ and the boundary $\Si$ are $C^2$ smooth, then there exists a constant $h_0>0$ such that for any $K\in\cM^\Ga\cup\cM^\Sigma$ satisfying $h_K< h_0$, we have $\eta_K\le Ch_K$ for some constant $C$ independent of $h_K$.
\end{lem}

\begin{proof} We only prove the lemma for the case of interface. The other case can be proved similarly. The argument extends that in Feistauer \cite[(3.27)]{F87} for the two-dimensional case. Since $\Ga$ is $C^2$ smooth, for any point $\bx\in\Ga$, there exist $r_{\bx}>0$ and a $C^2$ function $\Psi_{\bx}:\R^2\to\R$ such that, upon rotating and relabeling the coordinate axes if necessary, 
$\Om_1\cap B(\bx,r_{\bx})=\{\by\in B(\bx,r_{\bx}):y_3<\Psi_{\bx}(y_1,y_2)\}$, $\Om_2\cap B(\bx,r_{\bx})=\{\by\in B(\bx,r_{\bx}):y_3>\Psi_{\bx}(y_1,y_2)\}$, and $\Ga\cap B(\bx,r_{\bx})=\{\by\in B(\bx,r_{\bx}):y_3=\Psi_{\bx}(y_1,y_2)\}$.

By compactness, there exist finite number of points $\bx_j\in\Ga$, $1\le j\le \rev{M}$, such that the union of balls $B(\bx_j,r_{\bx_j}/2)$, $j=1,\cdots,\rev{M}$, covers the interface. Let $h_0=(\min_{1\le j\le \rev{M}}r_{\bx_j})/2$. Then any element $K\in\cM^\Ga$ satisfying $h_K<h_0$ is included in some ball $B(\bx_i,r_{\bx_i})$. Let $A,B,C$ be the three intersection points on the edges of $K$ that determines $\Ga^h_K$. If $\Pi_h\Psi_{\bx_i}:\R^2\to\R$ is the linear Lagrangian interpolation of $\Psi_{\bx_i}$ at the points $A,B,C$, then $\dist_{\rm H}(\Ga_K,\Ga_K^h)=\|\Psi_{\bx_i}-\Pi_h\Psi_{\bx_i}\|_{L^\infty(\R^2)}$. The lemma now follows from the standard finite element interpolation estimates as $\dist(A_K^i,\Ga_K^h)\ge Ch_K$.
\end{proof}

Lemma \ref{lem:3.0} implies that the following assumption is not very restrictive in practical applications.

\medskip
{\bf Assumption (H2)}
The interface and the boundary deviation $\eta_K\le 1/2$ for all $K\in\cM^{\Ga}\cup\cM^\Si$.
\medskip

We recall the following one-dimensional inverse domain estimate proved in \cite[Lemma 2.3]{Ch20}.

\begin{lem}\label{lem:3.1}
Let $I=(-1,1)$ and $I_\lam=(-\lam,\lam)$, $\lam>1$. We have
\ben
\|g\|_{L^2(I_\lam\backslash\bar I)}^2\le \frac 12\left[(\lam+\sqrt{\lam^2-1})^{2p+1}-1\right]\|g\|_{L^2(I)}^2\ \ \ \ \forall g\in Q_p(I_\lam),
\een
where $Q_p(I_\lam)$ is the set of polynomials of order $p\ge 1$ in $I_\lam$.
\end{lem}

We remark that the growing factor $(\lam+\sqrt{\lam^2-1})^{2p+1}$ in the bound is sharp, which is attained by the Chebyshev polynomials. It is well known that the Chebyshev polynomials $C_n(t)=\frac 12[(t+\sqrt{t^2-1})^n+(t-\sqrt{t^2-1})^n]$, $n\ge 0$, see DeVore and Lorentz \cite[P.76]{De93}.

For any integer $p\ge 1$ and any Lipschitz domain $D\subset \R^3$, we denote $Q_p(D)$ the set of polynomials of order $p$ in each variable in $D$. The following lemma can be proved as in \cite[Lemma 2.4]{Ch20} by using Lemma \ref{lem:3.1}. We omit the details.

\begin{lem}\label{lem:3.2}
Let $T$ be a tetrahedron with vertices $O=(0,0,0)^T$, $A=(a_1,a_2,a_3)^T$, $B=(b_1,b_2,0)^T$, and $C=(c_1,c_2,0)^T$, where $a_3>0$. For any $\de\in (0,a_3)$, let the tetrahedron $T_\de=\{\bx\in T:\dist(\bx,\De OBC)>\de\}$. Then, we have
	\ben
	\|v\|_{L^2(T)}\le\mathsf{T}\left(\frac{1+\delta a_3^{-1}}{1-\de a_3^{-1}}\right)^{3p+3/2}\|v\|_{L^2(T_\de)}\ \ \ \ \forall v\in Q_p(T),
	\een
where $\mathsf{T}(t)=t+\sqrt{t^2-1}\ \ \forall t\ge 1$.
\end{lem}

The following inverse trace estimate on curved domains is the main result in this subsection.

\begin{lem}\label{lem:3.4}
Let Assumptions (H1) and (H2) be satisfied. Then for any $K\in\cM^\Ga\cup\cM^\Si$, $K_i=K\cap\Om_i$, $i=1,2$,
	\ben
	\|v\|_{L^2(\pa K_i)}\le Cph_K^{-1/2}\mathsf{T}\left(\frac{1+3\eta_K}{1-\eta_K}\right)^{3p}\|v\|_{L^2(K_i)}\
	\ \ \ \forall v\in Q_p(K), 
	\een
where the constant $C$ is independent of $h_K, p$, and $\eta_K$.
\end{lem}

\begin{proof} We only prove the theorem when $K\in\cM^\Ga$. The other case $K\in\cM^\Si$ can be proved similarly.
For $i=1,2$, let $K_i^h$ be the polyhedron defined in Lemma \ref{lem:3.1new}. $K_i^h$ is divided into tetrahedrons, each of them with one vertex at $A_K^i$ as follows
\ben
\bar K^h_i=(\cup^{n}_{j=1}\bar T_{ij})\cup(\cup_{j=n+1}^{m_i}\bar T_{ij}),
\een
where $T_{ij}$, $j=1,\cdots,n$, has at most two vertices in $\mathcal{V}_K$ and $T_{ij}$, $j=n+1,\cdots,m_i$, has three vertices in $\mathcal{V}_K$. By Lemma \ref{lem:3.1new}, each tetrahedron $T_{ij}$, $j=1,\cdots,m_i$, is shape regular. Let $\de=\dist_{\rm H}(\Ga_K,\Ga_K^h)$ and $d_i=\dist(A_K^i,\Ga_K^h)$ so that $\de/d_i\le\eta_K$. Let $\Ga_K^{h\pm \de}$ be two planes parallel to $\Ga_K^h$ whose distance to $A_K^i$ is $d_i\pm\de$.

For $j=n+1,\cdots,m_i$, denote $T^*_{ij}$ the infinite cone with vertex $A_K^i$ that agrees with $T_{ij}$ in a neighborhood of $A_K^i$.  Denote 
$T_{ij}^{h\pm \de}$ the finite part of the cone $T^*_{ij}$ terminated by $\Ga_K^{h\pm\de}$. Then
$T_{ij}^{h-\de}\subset T_{ij}\cap K_i\subset T^{h+\de}_{ij}$.
By Lemma \ref{lem:3.2}, we have for $j=n+1,\cdots, m_i$,
\be\label{e2}
\|v\|_{L^2(T_{ij}^{h+\de})}&\le&\mathsf{T}\left(\frac{1+2\de(d_i+\de)^{-1}}{1-2\de(d_i+\de)^{-1}}\right)^{3p+3/2}\|v\|_{L^2(T_{ij}^{h-\de})}\nn\\
&\le&\mathsf{T}\left(\frac{1+3\eta_K}{1-\eta_K}\right)^{3p+3/2}\|v\|_{L^2(K_i\cap T_{ij})}\ \ \ \ \forall v\in Q_p(K).
\ee
By the $hp$-inverse trace inequality in Warburton and Hesthaven \cite{Wa03}, we have then for $j=n+1,\cdots,m_i$,
\be\label{l1}
\|v\|_{L^2(\pa T_{ij}^{h+\de}\cap\pa K_i)}&\le&\|v\|_{L^2(\pa T_{ij}^{h+\de})}\nn\\&\le&Cph_K^{-1/2}\|v\|_{L^2(T_{ij}^{h+\de})}\nn\\
&\le&Cph_K^{-1/2}\mathsf{T}\left(\frac{1+3\eta_K}{1-\eta_K}\right)^{3p+3/2}\|v\|_{L^2(K_i\cap T_{ij})}\ \ \forall v\in Q_p(K).
\ee

For each tetrahedron $T_{ij}$. $j=1,\cdots,n$, denote by $B,C,D$ the other three vertices other than $A_K^i$, see Fig.\ref{fig:3.3} (right). If $B\in\mathcal{V}_K$, the set of intersection points of $\Ga$ and the edges of $K$, $\Ga_K^{h-\de}$ intersects the edge $A_K^iB$ at $B'$. Let the half line originated from $A_K^i$ that is perpendicular to $\Ga_K^h$ intersect $\Ga_K^{h-\de}$ at $A'$ and intersect the plane parallel to $\Ga_K^h$ on which $B$ is located at $A$, then $|A_K^iA'|=d_i-\de$, $|A_K^iA|\le d_i+\de$. Thus  
$|A_K^iA'|/|A_K^i A|\ge(d_i-\de)/(d_i+\de)$ and consequently,
$|BB'|/|A_K^i B|=|AA'|/|A_K^iA|\le 2\de/(d_i+\de)$. Similarly, $|CC'|/|A_K^iC|\le 2\de/(d_i+\de)$ if $C\in\mathcal{V}_K$ and $\Ga_K^{h-\de}$ intersects the edge $A_K^iC$ at $C'$. Notice that by definition, $T_{ij}$ has at most two vertices in $\mathcal{V}_K$.

Now let $B''C''D''$ be the intersection of the plane parallel to $BCD$ whose distance to $BCD$ is $\de_1=\frac{2\de}{d_i+\de}\,h_1$, where $h_1=\dist(A_K^i,\De BCD)$. Then the
tetrahedron $\hat T_{ij}$ with vertices $A_K^i, B'',C'',D''$ is inside $T_{ij}\cap K_i$. Again by Lemma \ref{lem:3.2}, we have for $j=1,\cdots, n$, 
\be\label{e3}
\|v\|_{L^2(T_{ij})}&\le&\mathsf{T}\left(\frac{1+\de_1h_1^{-1}}{1-\de_1h_1^{-1}}\right)^{3p+3/2}\|v\|_{L^2(\hat T_{ij})}\nn\\
&\le&\mathsf{T}\left(\frac{1+3\eta_K}{1-\eta_K}\right)^{3p+3/2}\|v\|_{L^2(K_i\cap T_{ij})}\ \ \ \ \forall v\in Q_p(K).
\ee

Next let $\Ga^{h+\de}_K$ intersect the extended half lines of $A_K^iB, A_K^iC$ at $B''',C'''$, respectively, see Fig.\ref{fig:3.3} (right). Denote by $\widetilde{\Delta A_K^iBC}$ the curved triangle on the plane spanned by $A_K^i, B,C$ bounded by two edges $A_K^iB, A_K^iC$, and one curved edge $\widetilde{BC}$ that is the intersection $\Ga$ and the plane spanned by $A_K^i, B,C$. Then by Lemma \ref{lem:3.2}, \eqref{e3}, and the inverse trace inequality, we have for $j=1,\cdots, n$,
\be\label{l2}
\|v\|_{L^2(\widetilde{\Delta A_K^iBC})}&\le&\|v\|_{L^2(\Delta A_K^iB'''C''')}\nn\\
&\le&\mathsf{T}\left(\frac{1+2\de(d_i+\de)^{-1}}{1-2\de(d_i+\de)^{-1}}\right)^{3p+3/2}\|v\|_{L^2(\Delta A_K^i B'C')}\nn\\
&\le&\mathsf{T}\left(\frac{1+3\eta_K}{1-\eta_K}\right)^{3p+3/2}\|v\|_{L^2(\pa T_{ij})}\nn\\
&\le&C ph_K^{-1/2}\mathsf{T}\left(\frac{1+3\eta_K}{1-\eta_K}\right)^{3p+3/2}\|v\|_{L^2(K_i\cap T_{ij})}.
\ee
Let $F_{ij}^l$, $l=1,\cdots, L_{ij}$, $L_{ij}\le 3$, be the faces of $T_{ij}$ included in $\pa K$, $j=1,\cdots, n$, with the convention that $F_{ij}^1$ has two vertices on the edges of $K$. Denote by $\tilde F_{ij}^1$ the curved triangle on $\pa K$ that agrees with $F_{ij}^1$ in the neighborhood of $A_K^i$. Set $\tilde F_{ij}^l=F_{ij}^l$, $l=2,\cdots,L_{ij}$. Then by \eqref{l2} for $\tilde F_{ij}^1$ and the inverse trace inequality for $\tilde F_{ij}^l, l=2,\cdots,L_{ij}$, and \eqref{e3}, we obtain for $j=1,\cdots,n$,
\be\label{l3}
\sum_{l=1}^{L_{ij}}\|v\|_{L^2(\tilde F_{ij}^l)}\le Cph_K^{-1/2}\mathsf{T}\left(\frac{1+3\eta_K}{1-\eta_K}\right)^{3p+3/2}\|v\|_{L^2(K_i\cap T_{ij})}.
\ee

We recall the following trace inequality in Xiao et al. \cite[Lemma 3.1]{Wang} 
\beq\label{e4}
\|v\|_{L^2(\Ga_K)}\le C\|v\|_{L^2(K_i)}^{1/2}\|v\|_{H^1(K_i)}^{1/2}+\|v\|_{L^2(\pa K_i\backslash\bar\Ga_K)}\ \ \ \ \forall v\in H^1(K_i),\ \ i=1,2,
\eeq
where the constant $C$ is independent of $h_K$. Since on each face $F$ of $K$, $\pa K_i\cap F$ is the union of one curved triangle that either is $\tilde F_{ij}^1$ for $j=1,\cdots,n$, or is included in $T^{h+\de}_{ij}$ for $j=n+1,\cdots,m_i$, and at most two triangles that are faces of $T_{ij}, j=1,\cdots,m_i$, we have 
\ben
\pa K_i\backslash\bar\Ga_K\subset\left(\cup^n_{j=1}(\cup^{L_{ij}}_{l=1}\tilde F_{ij}^l)\right)\cup\left(\cup^{m_i}_{j=n+1}(\pa T_{ij}^{h+\de}\cap\pa K_i)\right),
\een
Thus by \eqref{l1} and \eqref{l3}, we have
\be\label{e5}
\|v\|_{L^2(\pa K_i)}&\le&C\|v\|_{L^2(K_i)}^{1/2}\|v\|_{H^1(K_i)}^{1/2}+\sum^n_{j=1}\sum^{L_{ij}}_{l=1}\|v\|_{L^2(\tilde F_{ij}^l)}+
\sum^{m_i}_{j=n+1}\|v\|_{L^2(\pa T_{ij}^{h+\de}\cap\pa K_i)}\nn\\
%&\le&C\|v\|_{L^2(K_i)}^{1/2}\|v\|_{H^1(K_i)}^{1/2}+Cph_K^{-1/2}\sum^{n}_{j=1}\|v\|_{L^2(T_{ij})}\nn\\
&\le&C\|v\|_{L^2(K_i)}^{1/2}\|v\|_{H^1(K_i)}^{1/2}+Cph_K^{-1/2}\mathsf{T}\left(\frac{1+3\eta_K}{1-\eta_K}\right)^{3p+3/2}\|v\|_{L^2(K_i)}.
\ee
Finally, since $K_i\subset(\cup^{n}_{j=1}T_{ij})\cup(\cup^{m_i}_{j=n+1}T_{ij}^{h+\de})$, we obtain by the $hp$-inverse estimate (see, e.g., Schwab \cite[Theorem 4.76]{Sc98}), \eqref{e2}, and \eqref{e3} that
\be\label{e6}
\|\na v\|_{L^2(K_i)}&\le&Cp^2h_K^{-1}\left(\sum^{n}_{j=1}\|v\|_{L^2(T_{ij})}+\sum^{m_i}_{j=n+1}\|v\|_{L^2(T_{ij}^{h+\de})}\right)\nn\\
&\le&Cp^2h_K^{-1}\mathsf{T}\left(\frac{1+3\eta_K}{1-\eta_K}\right)^{3p+3/2}\|v\|_{L^2(K_i)}.
\ee
This completes the proof by using \eqref{e5}.
\end{proof}

%%%%%%%%%%%%%%%%%%%%%%%%%%%%%%%%%%%%%%%%%%%%%%%%%
\subsection{The unfitted finite element method}
%%%%%%%%%%%%%%%%%%%%%%%%%%%%%%%%%%%%%%%%%%%%%%%%%%

We introduce the finite element space using the idea of ``doubling of unknowns" in Hansbo and Hansbo \cite{Ha02}.  For any integer $p\ge 1$, we define the unfitted finite element spaces 
\ben
& &\bX_p(\cM)=\{\bv=\bv_1\chi_{\Om_1}+\bv_2\chi_{\Om_2}:\bv_i|_K\in \bQ_p(K), i=1,2\},\\
& &M_p(\cM)=\{q=q_1\chi_{\Om_1}+q_2\chi_{\Om_2}:q_i\in Q_p(K),i=1,2\}.
\een
Clearly, $\bX_p(\cM)=M_p(\cM)^3$. \rev{We remark that one may use $P_p(K)$ instead of $Q_p(K)$ in the definition of the unfitted finite element spaces. We use $Q_p(K)$ because this choice offers the possibility to construct conforming finite elements away from the interface or the boundary.}

Let $\cF^{\rm int}$ denote the set of interior faces of the mesh $\cM$, $\cF^\Ga=\cup_{K\in\cM}\Ga_K$, and $\cF^{\Sigma}=\cup_{K\in\cM}\Si_K$. Since hanging nodes are allowed, $F\in\cF^{\rm int}$ can be part of a face of an adjacent element. We set $\cF=\cF^{\rm int}\cup\cF^\Ga$ and $\overline{\cF}= \cF^{\rm int}\cup\cF^\Ga\cup\cF^{\Sigma}$.
For any subset $\widehat\cM\subset\cM$ and $\widehat\cF\subset\overline{\cF}$, we use the notation
\ben
(u,v)_{\widehat\cM}:=\sum_{K\in\widehat\cM}(u,v)_K,\ \ \la u,v\ra_{\widehat\cF}:=\sum_{F\subset\widehat\cF}\la u,v\ra_F,
\een
where $(u,v)_K$ is the inner product of $L^2(K)$ and $\la u,v\ra_F$ is the inner product of $L^2(F)$.

For any $F\in\overline{\cF}$, we fix a unit normal vector $\bn_F$ of $F$ with the convention that $\bn_F$ is the unit outer normal to $\Sigma$ if $F\in\cF^\Sigma$ and $\bn_F$ is the unit outer normal to $\pa\Om_1$ if $F\in\cF^\Ga$. We define the piecewise constant normal vector function $\bn\in L^\infty(\overline{\cF})^3$ by
$\bn|_F=\bn_F\ \ \forall F\in\overline{\cF}$. 

For any $v\in H^1(\cM):=\{v=v_1\chi_{\Om_1}+v_2\chi_{\Om_2}:v_i|_K\in H^1(K), i=1,2\}$, we define the jump of $v$ across $F$ as 
\ben
\lj v\rj_F:=v\rev{^-}-v\rev{^+}\ \ \forall F\in\cF,\ \ \ \
\lj v\rj_F:=v\rev{^-}\ \ \forall F\in\cF^\Sigma,
\een
where $v\rev{^\pm}$ is the trace of $v$ on $F$ in the $\pm \bn_F$ direction.

Denote $\bW=\mu^{-1}\na\times\bE$. Then \eqref{a1}-\eqref{a2} can be rewritten as
\ben
& &\mu\bW-\na\times\bE=0,\ \ \na\times\bW-k^2\vep\bE-\vep\na\vp=\bJ,\ \ \div(\vep\bE)=0\ \ \ \ \mbox{in }\Om,\\
& &\bE\times\bn=\bg\times\bn,\ \ \vp=0\ \ \ \ \mbox{on }\Sigma.
\een
Following Cockburn and Shu \cite{Co98}, Alvarado and Castillo \cite{Al16}, the LDG method is to find $(\bW_h,\bE_h,\vp_h)\in\bX_p(\cM)\times\bX_p(\cM)\times M_p(\cM)$ such that for any $K\in\cM$,
\be
& &(\mu\bW_h,\bt)_K-(\bE_h,\na\times\bt)_K-\la\bn_K\times\rev{\widehat{\bE}_h},\bt\ra_{\pa K}
-\la\bn_K\times\rev{\widehat{\bE}_h},\lj\bt\rj\ra_{\Ga_K}=0,\label{f1}\\
& &(\bW_h,\na\times\bv)_K+\la\bn_K\times\rev{\ddhat{\bW}_h},\bv_T\ra_{\pa K}+\la\bn_K\times\rev{\ddhat{\bW}_h},\lj\bv_T\rj\ra_{\Ga_K}-k^2(\vep\bE_h,\bv)_K\nn\\
& &\hskip1cm +(\vp_h,\div(\vep\bv))_K-\la\rev{\widehat{\vp}_h},\vep\bv\cdot\bn_K\ra_{\pa K}-\la\rev{\widehat{\vp}}_h,\lj\vep\bv\cdot\bn\rj\ra_{\Ga_K}=(\bJ,\bv)_K,\label{f2}\\
& &(\vep\bE_h,\na q)_K-\la\rev{\ddhat{\vep\bE}_h}\cdot\bn_K,q\ra_{\pa K}-\la\rev{\ddhat{\vep\bE}_h}\cdot\bn,\lj q\rj\ra_{\Ga_K}=0,\label{f3}
\ee
for all $(\bt,\bv,q)\in\bX_p(\cM)\times\bX_p(\cM)\times M_p(\cM)$. \rev{Here we define the numerical fluxes
\ben
\widehat{\bE}_h|_F=\bE^+_h,\ \widehat{\vp}_h|_F=\vp^+_h,\ \ddhat{\bW}_h|_F=\bW_h^-,\ \ddhat{\vep\bE}_h|_F=(\vep\bE_h)^-\ \ \ \ \mbox{on } F\in\cF,
\een}
and due to the boundary conditions $\bn\times\bE=\bn\times\bg$, $\vp=0$ on $\Sigma$,
\beq\label{f4}
\bn\times\rev{\widehat{\bE}_h}|_F=\bn\times\bg,\ \ \rev{\widehat{\vp}_h}|_F=0\ \ \ \ \mbox{on } F\in\cF^{\Sigma}.
\eeq
Now integrating by parts in \eqref{f1}, summing the equations over $K\in\cM$, and using the following elementary DG magic formula
that for any $a,b\in H^1(\cM)$, 
\ben
\lj ab\rj_F=a^-\lj b\rj_F+\lj a\rj_Fb^+=a^+\lj b\rj_F+\lj a\rj_Fb^-\ \ \ \ \forall F\in\cF,
\een
one can obtain
\beq
(\mu\bW_h,\bt)_\cM-(\na_h\times\bE_h,\bt)_\cM+\la\lj\bn\times\bE_h\rj,\bt^-\,\ra_{\overline{\cF}}-\la\bn\times\bg,\bt^-\,\ra_{\cF^\Sigma}=0,\label{f6}
\eeq
where $\na_h\times\bE_h|_K=(\na\times\bE_h|_{K_1})\chi_{\Om_1}+(\na\times\bE_h|_{K_2})\chi_{\Om_2}$ on each element $K\in\cM$. Similarly, from \eqref{f2} we deduce
\be\label{f7}
& &(\bW_h,\na_h\times\bv)_\cM+\la\bn\times\bW_h^-,\lj\bv_T\rj\ra_{\overline{\cF}}-k^2(\vep\bE_h,\bv)_\cM\nn\\
& &\hskip2cm -(\na_h\vp_h,\vep\bv)_\cM+\la\lj\vp_h\rj,(\vep\bv)^-\cdot\bn\ra_{\overline{\cF}}=(\bJ,\bv)_\cM.
\ee
From \eqref{f3} we obtain
\beq\label{f8}
(\vep\bE_h,\na_h q)_\cM-\la(\vep\bE_h)^-\cdot\bn,\lj q\rj\ra_{\overline{\cF}}=0.
\eeq
Define the lifting operators $\mathsf{L}:\bH^1(\cM)\to\bX_p(\cM)$, $\mathsf{L}_1:\bL^2(\Sigma)\to\bX_p(\cM)$ by
\beq\label{m2}
(\mathsf{L}(\bv),\bt)_\cM=\la\lj\bn\times\bv\rj,\bt^-\,\ra_{\overline{\cF}},\ \ (\mathsf{L}_1(\bg),\bt)_\cM=\la\bn\times\bg,\bt^-\,\ra_{\cF^\Sigma}\ \ \ \ \forall\bt\in\bX_p(\cM).
\eeq
Then \eqref{f6} yields $\mu\bW_h=\na_h\times\bE_h-\mathsf{L}(\bE_h)+\mathsf{L}_1(\bg)$ and thus from \eqref{f7}, for any $\bv\in\bX_p(\cM)$,
\be\label{f9}
& &(\mu^{-1}(\na_h\times\bE_h-\mathsf{L}(\bE_h)),\na_h\times\bv-\mathsf{L}(\bv))_\cM-k^2(\vep\bE_h,\bv)_\cM-(\na_h\vp_h,\vep\bv)_\cM\nn\\
& &\hskip1cm +\la\lj\vp_h\rj,(\vep\bv)^-\cdot\bn\ra_{\overline{\cF}}=(\bJ,\bv)_\cM
-(\mu^{-1}\mathsf{L}_1(\bg),\na_h\times\bv-\mathsf{L}(\bv))_\cM.
\ee
In practical computations, \eqref{f8} and \eqref{f9} must be implemented by adding stablization or penalty terms.

Introduce the sesquilinear forms $a:\bH^1(\cM)\times\bH^1(\cM)\to\C$, $b:\bH^1(\cM)\times H^1(\cM)\to\C$, and $c:H^1(\cM)\times H^1(\cM)\to\C$ as follows:
\ben
a(\bu,\bv)&=&(\mu^{-1}(\na_h\times\bu-\mathsf{L}(\bu)),\na_h\times\bv-\mathsf{L}(\bv))_\cM\\
& &\,+\,\la\al\tau^{-1}\lj\bn\times\bu\rj,\lj\bn\times\bv\rj\ra_{\overline{\cF}}+\la\al\tau^{-1}\lj\vep\bu\cdot\bn\rj,\lj\vep\bv\cdot\bn\rj\ra_{{\cF}}\\
b(\bv,q)&=&(\vep\bv,\na_h q)_\cM-\la(\vep\bv)^-\cdot\bn,\lj q\rj\ra_{\overline{\cF}},\\
c(\vp,q)&=&(\tau\tau_1\na_h\vp,\na_h q)_\cM+\la\tau\lj\vp\rj,\lj q\rj\ra_{\overline{\cF}}.
\een
Here for any $F\in\overline{\cF}$, $\al|_F=\al_0\Theta_F$, $\tau|_F=h_Fp^{-2}$, and $\tau_1|_F=h_Fp^{-1}$, where $\al_0>0$ is some fixed constant, 
\ben
\Theta_F=\max_{F\cap\bar K\not=\emptyset, K\in\cM}\,\Theta_K,\ \ \Theta_K=\left\{\begin{array}{ll}
\mathsf{T}\left(\frac{1+3\eta_K}{1-\eta_K}\right)^{6p} & \mbox{if }K\in\cM^\Ga\cup\cM^\Sigma,\cr\\
1 & \mbox{otherwise},
\end{array}\right.
\een
and $h_F=(h_K+h_{K'})/2$ if $F=\pa K\cap\pa K'\in\cF^{\rm int}$ for some elements $K,K'\in\cM$ and $h_F=h_K$ if $F=\Ga_K\in\cF^\Ga$ or $F=\Sigma_K\in\cF^\Sigma$ for some $K\in\cM$. 

The penalty is to penalize $\lj\bn\times\bu\rj$ and $\lj\vep\bu\cdot\bn\rj$ in the $H^{1/2}$ norm. The Lagrangian multiplier $\vp$ is penalized in the $L^2$ norm so that $\na_h\vp$ in the $H^{-1}$ and $\lj\vp\rj$ in the $H^{-1/2}$ norm. Lemma \ref{lem:4.3} below shows that the penalty of $\lj\vep\bu\cdot\bn\rj$ in the $H^{1/2}$ norm
implies the penalty of $\div(\vep\bu)$ in the $L^2$ norm, which we have not included explicitly in our sesquilinear form $a(\cdot,\cdot)$ to reduce the computational costs. 

\revto{We remark that the penalty of $\lj\bn\times\bu\rj$ in the $H^{1/2}$ norm is widely used in the literature (see, e.g., Houston et al. \cite{Ho04}, \cite{Ho05}, and Bonito et al. \cite{Bo16}) in which conforming subspaces of the broken space are used to derive optimal convergence rates. However, since the trace space of $\bH(\curl;\Om_1\cup\Om_2)$ on the interface $\Ga$ is $\bH^{-1/2}(\div_{\Ga};\Ga)$, the natural penalty for $\lj\bn\times\bu\rj$ would be the $H^{-1/2}$ norm, which leads to optimal convergence even for solutions with lower regularity, see Brenner et al. \cite{Brenner09} for a nonconforming penalty method, Beur\~{a}o da Veiga et al. \cite{Bei22} and Cao et al. \cite{Cao21} for virtual element methods. The penalties for $\lj\vep\bu\cdot\bn\rj$, $\lj\vp\rj$, $\na_h\vp$ used in this paper are a special case of that in \cite{Bo16}, in which an interior penalty method with $C^0$ finite elements using the penalty of $\div(\vep\bu)$ in the $H^{-\gamma}$, $\vp$ in the $H^{\gamma}$, and $\lj\vep\bu\cdot\bn\rj$ in the $H^{-\gamma+1/2}$ norm, where $0\le\gamma\le 1$, is studied. Another different penalty
of the divergence free condition for Maxwell equations is considered in Costabel and Dauge \cite{Co02}. We also remark that unfitted finite element methods for Maxwell interface problems with piecewise smooth interfaces, whose solutions have lower regularities, require further investigations.}

It is clear that by integrating by parts and using the DG magic formula, we obtain
\beq\label{xx2}
b(\bv,q)=-(\na_h\cdot(\vep\bv),q)_\cM+\la\lj\vep\bv\cdot\bn\rj,q^+\ra_\cF\ \ \ \ \forall (\bv,q)\in\bH^1(\cM)\times H^1(\cM).
\eeq

The unfitted finite element method for solving \eqref{a1}-\eqref{a2} is to find $(\bE_h,\vp_h)$ $\in\bX_p(\cM)\times M_p(\cM)$ such that
\be
& &a(\bE_h,\bv)-\overline{b(\bv,\vp_h)}-k^2(\vep\bE_h,\bv)_\cM=\bF_h(\bv)\ \ \forall\bv\in\bX_p(\cM),\ \ \label{g1}\\
& &b(\bE_h,q)+c(\vp_h,q)=0\ \ \forall q\in M_p(\cM),\label{g2}
\ee
where $\bF_h(\bv)=(\bJ,\bv)_\cM-(\mu^{-1}\mathsf{L}_1(\bg),\na_h\times\bv-\mathsf{L}(\bv))_\cM+\la\al\tau^{-1}\bn\times\bg,\lj\bn\times\bv\rj\ra_{\cF^\Sigma}\ \ \forall\bv\in\bX_p(\cM)$.
The well-posedness of the problem and the convergence of $(\bE_h,\vp_h)$ to the solution $(\bE,\vp)$ of \eqref{a1}-\eqref{a2} will be considered in \S4.2.

%%%%%%%%%%%%%%%%%%%%%%%%%%%%%%%%%%%%%%%%%%%%%%%%%%%
\section{The finite element convergence analysis}
%%%%%%%%%%%%%%%%%%%%%%%%%%%%%%%%%%%%%%%%%%%%%%%%%%%

In this section, we study the well-posedness and the convergence of the unfitted finite element method for the time-harmonic Maxwell equations in \eqref{g1}-\eqref{g2}. We first prove optimal error estimates for a projection operator for the coercive Maxwell equations in \S4.1, which is used to prove the optimal error estimates for our unfitted finite element method using the Schatz argument in \S4.2. 

In this section, unless otherwise stated, we will denote $C>0$ the generic constant which may depend on the wave number $k$ and the coefficients $\vep_i,\mu_i$, $i=1,2$, but is independent of $p$, $h_K$ for all $K\in\cM$, and $\eta_K$ for all $K\in\cM^\Ga\cup\cM^\Si$.

%%%%%%%%%%%%%%%%%%%%%%%%%%%%%%%%%%%%%%%%%%%%%%%%%%%
\subsection{The projection operator for coercive Maxwell equations}
%%%%%%%%%%%%%%%%%%%%%%%%%%%%%%%%%%%%%%%%%%%%%%%%%%%

In this subsection we study a projection operator for the $\bH(\curl)$-coercive Maxwell equations\rev{,} which will be used
in our analysis of the problem \eqref{g1}-\eqref{g2} in the next subsection. We first introduce the following DG norms: for any $\bv\in\bH^1(\cM)$, $q\in H^1(\cM)$,
\ben
& &|\bv|_{\bX_p(\cM)}^2:=\|\na_h\times\bv\|_\cM^2+\|\al^{1/2}\tau^{-1/2}\lj\bn\times\bv\rj\|_{\overline{\cF}}^2
+\,\|\al^{1/2}\tau^{-1/2}\lj\vep\bv\cdot\bn\rj\|_{{\cF}}^2,\\
& &\|\bv\|_{\bX_p(\cM)}^2:=|\bv|^2_{\bX_p(\cM)}+\|k\bv\|_\cM^2,\\
& &\|q\|_{M_p(\cM)}^2:=\|\tau^{1/2}\tau_1^{1/2}\na_h q\|_\cM^2+\|\tau^{1/2}\lj q\rj\|_{\overline{\cF}}^2.
\een
We define the projection operator $\Pi_h:\bH^1(\cM)\times H^1_0(\Om)\to\bX_p(\cM)\times M_p(\cM)$ as follows:  
for any $(\bu,\phi)\in\bH^1(\cM)\times H^1_0(\Om)$, $\Pi_h(\bu,\phi)=(\bu_h,\phi_h)\in\bX_p(\cM)\times M_p(\cM)$ satisfies
\be
& &a(\bu_h-\bu,\bv)+k^2(\vep(\bu_h-\bu),\bv)_\cM-\overline{b(\bv,\phi_h-\phi)}=0\ \ \forall\bv\in\bX_p(\cM),\ \ \ \ \label{g3}\\
& &b(\bu_h-\bu,q)+c(\phi_h-\phi,q)=0\ \ \ \ \forall q\in M_p(\cM).\label{g4}
\ee

It is easy to see that the operator $\Pi_h$ is well defined. In fact, we only need to show the uniqueness, for that purpose, we let $(\bu,\phi)=(0,0)$. By taking $\bv=\bu_h$ in \eqref{g3}, $q=\phi_h$ in \eqref{g4}, and then adding the first equation with the complex conjugate of the second equation, we obtain easily $\bu_h=0, \phi_h=0$. This shows the uniqueness of the solution to \eqref{g3}-\eqref{g4}. \rev{We remark that thanks to the stabilization term $c(\cdot,\cdot)$, the stability of \eqref{g3}-\eqref{g4} can be easily proved without considering the inf-sup condition of $b(\cdot,\cdot)$.}

The main purpose of this subsection is to show an $hp$-error estimate for 
$(\bu,\phi)-\Pi_h(\bu,\phi)$. We first recall the following $hp$-approximation result in Babu\v{s}ka and Suri \cite[Lemma 4.5]{Ba87}, Melenk \cite[Lemma B.3]{Me05}. 

\begin{lem}\label{lem:4.1}
Let $s\ge 0$ and $p\ge 1$. For any $K\in\cM$, there exists an interpolation operator $\pi_{hp}^K:H^s(K)\to Q_p(K)$ such that for any $u\in H^s(K)$, 
\ben
\|u-\pi^{hp}_K(u)\|_{H^j(K)}\le C\frac{h_K^{\nu-j}}{p^{s-j}}\|u\|_{H^s(K)},\ \ 0\le j\le s,
\een
where $\nu=\min(p+1,s)$ and the constant $C$ is independent of $h_K,p$, but may depend on $s$.
\end{lem}

By the multiplicative trace inequality
\beq\label{m1}
\|v\|_{\rev{L}^2(\pa K)}\le Ch_K^{-1/2}\|v\|_{L^2(K)}+C\|v\|_{L^2(K)}^{1/2}\|v\|_{H^1(K)}^{1/2}\ \ \ \ \forall v\in H^1(K),
\eeq
and the inequality \eqref{e4}, we obtain from Lemma \ref{lem:4.1} that for $i=1,2$,
\be\label{g5}
& &\|u-\pi^{hp}_K(u)\|_{L^2(\pa K_i)}\nn\\
&\le&C(\|u-\pi^{hp}_K(u)\|_{L^2(K)}^{1/2}\|\na(u-\pi^{hp}_K(u))\|_{L^2(K)}^{1/2}+\|u-\pi^{hp}_K(u)\|_{L^2(\pa K)})\nn\\
&\le&C\frac{h_K^{\nu-1/2}}{p^{s-1/2}}\|u\|_{H^s(K)},\ \ s\ge 1.
\ee

\begin{lem}\label{lem:4.2}
There exists a constant $c_{\rm coer}>0$ such that $|a(\bv,\bv)|\ge c_{\rm coer}|\bv|_{\bX_p(\cM)}^2$ for all $\bv\in\bX_p(\cM)$. The constant $c_{\rm coer}$
is independent of $p$, $h_K$ for all $K\in\cM$, and $\eta_K$ for all $K\in\cM^\Ga\cup\cM^\Si$.
\end{lem}

\begin{proof} By the inverse trace estimate in Lemma \ref{lem:3.4}, there exists a constant $c_L>0$ independent of $p$, $h_K$ for all $K\in\cM$, and $\eta_K$ for all $K\in\cM^\Ga\cup\cM^\Si$ such that 
\beq\label{m3}
\|\mathsf{L}(\bu)\|_\cM\le 
c_L\|\al^{1/2}\tau^{-1/2}\lj\bn\times\bu\rj\|_{\overline{\cF}}\ \ \ \ \forall \bu\in\bH^1(\cM). 
\eeq
The rest of the proof follows easily from the standard argument (see, e.g., \cite[Theorem 2.1]{Ch20}) with $c_{\rm coer}=(4+c_L^2)^{-1}$.
\end{proof}

\begin{lem}\label{lem:4.5}
Let $s\ge 1$. There exists an interpolation operator $I_{h,p}:H^s(\Om_1\cup\Om_2)\to M_p(\cM)$ such that for any $u\in H^s(\Om_1\cup\Om_2)$, $i=1,2$,
\ben
& &\|u-I_{h,p}(u)\|_{H^j(K_i)}\le C\frac{h_K^{\nu-j}}{p^{s-j}}\|\tilde u_i\|_{H^s(K)},\ \ 0\le j\le s,\\
& &\|u-I_{h,p}(u)\|_{L^2(\pa K_i)}\le C\frac{h_K^{\nu-1/2}}{p^{s-1/2}}\|\tilde u_i\|_{H^s(K)},
\een
where $\tilde u_i\in H^s(\R^3)$ is the extension of $u|_{\Om_i}$ such that $\|\tilde u_i\|_{H^s(\R^3)}\le C\|u\|_{H^s(\Om_i)}$.
\end{lem}

\begin{proof} For $u=u_1\chi_{\Om_1}+u_2\chi_{\Om_2}$, $u_i\in H^s(\Om_i)$, $i=1,2$, let $\tilde u_i\in H^s(\R^3)$ be the Stein extension (see, e.g., Adams and Fournier \cite[Theorem 5.24]{Ad09}) of $u_i$ such that $\|\tilde u_i\|_{H^s(\R^3)}\le C\| u_i\|_{H^s(\Om_i)}$. We define $I_{h,p}(u)|_K=\pi_{hp}^K(\tilde u_1)\chi_{\Om_1}+\pi_{hp}^K(\tilde u_2)\chi_{\Om_2}\ \ \forall K\in\cM$, where $\pi_{hp}^K(\tilde u_i)$ is defined in Lemma \ref{lem:4.1}. The lemma follows easily by using Lemma \ref{lem:4.1} and \eqref{g5}.
\end{proof}

The following theorem is the main result of this subsection.

\begin{thm}\label{thm:4.1} 
Let $(\bu,\phi)\in\bH^{s}(\Om_1\cup\Om_2)\times H^{s}(\Om_1\cup\Om_2)$, $s\ge 1$, we have 
\ben
& &\|\bu-\bu_h\|_{\bX_p(\cM)}+\|\phi-\phi_h\|_{M_p(\cM)}\\
&\le&C\revt{\max_{K\in\cM}\left(\frac{\Theta_K^{1/2}h_K^{\nu-1}}{p^{s-3/2}}\right)}\left(\|\bu\|_{H^s(\Om_1\cup\Om_2)}+\|\phi\|_{H^s(\Om_1\cup\Om_2)}\right).
\een
\end{thm}

\begin{proof} Let $(\bu_I,\phi_I)=(I_{h,p}(\bu),I_{h,p}(\phi))\in\bX_p(\cM)\times M_p(\cM)$. By Lemma \ref{lem:4.5} 
we obtain easily that
\be
& &\|\bu-\bu_I\|_{\bX_p(\cM)}\le C\revt{\max_{K\in\cM}\left(\frac{\Theta_K^{1/2}h_K^{\nu-1}}{p^{s-3/2}}\right)}\|\bu\|_{H^s(\Om_1\cup\Om_2)},\label{hh1}\\
& &\|\phi-\phi_I\|_{M_p(\cM)}\le C\frac{h^{\nu}}{p^{s+1/2}}\|\phi\|_{H^s(\Om_1\cup\Om_2)},\label{hh2}
\ee
where $h=\max_{K\in\cM}h_K$.
From \eqref{g3}-\eqref{g4} we know that
\be
& &a(\bu_h-\bu_I,\bv)+k^2(\vep(\bu_h-\bu_I),\bv)-\overline{b(\bv,\phi_h-\phi_I)}\nn\\
& &\hskip1cm=a(\bu-\bu_I,\bv)+k^2(\vep(\bu-\bu_I),\bv)
-\overline{b(\bv,\phi-\phi_I)}\ \ \ \ \forall\bv\in\bX_p(\cM),\label{h2}\\
& &b(\bu_h-\bu_I,q)+c(\phi_h-\phi_I,q)=b(\bu-\bu_I,q)+c(\phi-\phi_I,q)\ \ \ \ \forall q\in M_p(\cM),\label{h3}
\ee
By \eqref{m3} and the definition of the DG norms, we have 
\ben
& &|a(\bu-\bu_I,\bv)+k^2(\vep(\bu-\bu_I),\bv)_\cM|\le C\|\bu-\bu_I\|_{\bX_p(\cM)}\|\bv\|_{\bX_p(\cM)},\\
& &|c(\phi-\phi_I,q)|\le \|\phi-\phi_I\|_{M_p(\cM)}\|q\|_{M_p(\cM)}.
\een
Next, by Lemma \ref{lem:3.4} and the fact that $\phi=0$ on $\Sigma$, it is easy to see that 
\ben
& &|b(\bv,\phi-\phi_I)|\le C\|\bv\|_\cM(\|\na_h(\phi-\phi_I)\|_\cM+\|\al^{1/2}\tau^{-1/2}\lj\phi-\phi_I\rj\|_{\overline{\cF}}),\\
& &|b(\bu-\bu_I,q)|\le C(\|\tau^{-1/2}\tau_1^{-1/2}(\bu-\bu_I)\|_\cM+\|\tau^{-1/2}(\vep\bu-\vep\bu_I)^-\cdot\bn\|_{\overline{\cF}})\|q\|_{M_p(\cM)}.
\een
Now by taking $\bv=\bu_h-\bu_I$ in \eqref{h2}, $q=\phi_h-\phi_I$ in \eqref{h3}, and using the standard argument, we have
\ben
& &\|\bu_h-\bu_I\|_{\bX_p(\cM)}+\|\phi_h-\phi_I\|_{M_p(\cM)}\\
&\le&C(\|\bu-\bu_I\|_{\bX_p(\cM)}+\|\tau^{-1/2}\tau_1^{-1/2}(\bu-\bu_I)\|_\cM+\|\tau^{-1/2}(\bu-\bu_I)\|_{\overline{\cF}})\\
& &\ +\,C(\|\phi-\phi_I\|_{M_p(\cM)}+\|\na_h(\phi-\phi_I)\|_\cM+\|\al^{1/2}\tau^{-1/2}\lj\phi-\phi_I\rj\|_{\overline{\cF}})\\
&\le&C\revt{\max_{K\in\cM}\left(\frac{\Theta_K^{1/2}h_K^{\nu-1}}{p^{s-3/2}}\right)}\left(\|\bu\|_{H^s(\Om_1\cup\Om_2)}+\|\phi\|_{H^s(\Om_1\cup\Om_2)}\right),
\een
where we have used Lemma \ref{lem:4.5} and \eqref{hh1}-\eqref{hh2}. This completes the proof by using \eqref{hh1}-\eqref{hh2}.
\end{proof}

We remark that the error estimate is optimal in $h$ and slightly suboptimal in $p$ under the regularity assumption $\bu\in\bH^s(\Om_1\cup\Om_2)$. For the $\bH(\curl)$ conforming finite element methods on conforming meshes for solving Maxwell equations with smooth coefficients, the optimal error estimate of the order $h^{\min(p,s)}/p^{s}$ can be obtained by using the projection-based interpolation operators in Demkowicz \cite{De08}, Melenk and Rojik \cite{Me19} under the regularity assumption $\bu\in\bH^s(\curl;\Om)=\{\bu\in\bH^s(\Om):\na\times\bu\in\bH^s(\Om)\}$.

%%%%%%%%%%%%%%%%%%%%%%%%%%%%%%%%%%%%%%%%%%%%%%%%%%%
\subsection{The time-harmonic Maxwell equations}
%%%%%%%%%%%%%%%%%%%%%%%%%%%%%%%%%%%%%%%%%%%%%%%%%%%

We start by studying the consistency errors of the unfitted finite element method.

\begin{lem}\label{lem:5.1}
Let $\bE\in\bH^s(\Om_1\cup\Om_2)$, $s\ge 2$, be the solution of the problem \eqref{a1}-\eqref{a2} and $(\bE_h,\vp_h)\in\bX_p(\cM)\times M_p(\cM)$ be the solution of \eqref{g1}-\eqref{g2}. Then 
\be
& &a(\bE-\bE_h,\bv)+\overline{b(\bv,\vp_h)}-k^2(\vep(\bE-\bE_h),\bv)_\cM=R(\bE,\bv)\ \ \forall\bv\in\bX_p(\cM),\label{x2}\\
& &b(\bE-\bE_h,q)-c(\vp_h,q)=0\ \ \ \ \forall q\in M_p(\cM),\label{x3}
\ee
where for any $\bv\in\bX_p(\cM)$, the residual $R(\bE,\bv)$ satisfies
\beq\label{x1}
|R(\bE,\bv)|\le C\revt{\max_{K\in\cM}\left(\frac{\Theta_K^{1/2}h_K^{\nu-1}}{p^{s-1}}\right)}\|\bE\|_{H^s(\Om_1\cup\Om_2)}\|\tau^{-1/2}\lj\bn\times\bv\rj\|_{\overline{\cF}}.
\eeq
\end{lem}

\begin{proof} The proof is similar to that in \S3.2 to derive the DG method. We multiply the first equation in \eqref{p1} by $\bv\in\bX_p(\cM)$ and integrate on each element $K\in\cM$,
\ben
& &(\mu^{-1}\na\times\bE,\na_h\times\bv)_K+\la\bn_K\times(\mu^{-1}\na\times\bE),\bv\ra_{\pa K}\\
& &\hskip1cm+\,\la\bn\times(\mu^{-1}\na\times\bE),\lj\bv\rj\ra_{\Ga_K}-k^2(\vep\bE,\bv)_K=(\bJ,v)_K.
\een
Summing the equations over all elements we obtain
\ben
(\mu^{-1}\na\times\bE,\na_h\times\bv-\mathsf{L}(\bv))_\cM-k^2(\vep\bE,\bv)
=(\bJ,\bv)+R(\bE,\bv),
\een
where $R(\bE,\bv)$ is defined as
\ben
R(\bE,\bv)=\la\mu^{-1}\na\times\bE,\lj\bn\times\bv\rj\ra_{\overline{\cF}}-(\mu^{-1}\na\times\bE,\mathsf{L}(\bv))_\cM\ \ \ \ \forall\bv\in\bX_p(\cM).
\een
From the definition \eqref{m2}, we know that $\mathsf{L}(\bE)=\mathsf{L}_1(\bg)$. This shows \eqref{x2} by \eqref{g1}.

Next we multiply the second equation in \eqref{p1} by $q\in M_p(\cM)$ and 
integrate on each element $K\in\cM$,
\ben
(\vep\bE,\na q)_K-\la\vep\bE\cdot\bn_K,q\ra_{\pa K}-\la\vep\bE\cdot\bn,\lj q\rj\ra_{\Ga_K}=0.
\een
Summing the equations over all elements yields $b(\bE,q)=0$ for any $q\in M_p(\cM)$. This shows \eqref{x3} by using \eqref{g2}.

It remains to show the estimate \eqref{x1} for the residual. Denote $\bW=\mu^{-1}\na\times\bE\in\bH^{s-1}(\Om_1\cup\Om_2)$ and $\bW_I=I_{h,p}(\bW)\in\bX_{p}(\cM)$ the interpolation function defined in Lemma \ref{lem:4.5}. By the definition of the lifting operator $\mathsf{L}$ in \eqref{m2} we know that
\ben
|R(\bE,\bv)|&=&|\la\bW-\bW_I^-,\lj\bn\times\bv\rj\ra_{\overline{\cF}}-(\bW-\bW_I,\mathsf{L}(\bv))_\cM|.
\een
\revt{Let $\sigma\in L^\infty(\Om)$ be a piecewise constant function $\sigma|_K=\Theta_K\ \forall K\in\cM$. By taking $\bt=\sigma^{-1}\mathsf{L}(\bv)\in\bX_p(\cM)$ in the definition of the lifting operator $\mathsf{L}$ in \eqref{m2}, one obtains easily by Lemma \ref{lem:3.4} that $\|\sigma^{-1/2}\mathsf{L}(\bv)\|_\cM\le C\|\tau^{-1/2}\lj\bn\times\bv\rj\|_{\overline{\cF}}$. Thus
\ben
|R(\bE,\bv)|\le C(\|\tau^{1/2}(\bW-\bW_I^-)\|_{\overline{\cF}}+\|\sigma^{1/2}(\bW-\bW_I)\|_\cM)\|\tau^{-1/2}\lj\bn\times\bv\rj\|_{\overline{\cF}}.
\een}
This shows \eqref{x1} by Lemma \ref{lem:4.5} and thus completes the proof.
\end{proof}

\begin{lem}\label{lem:4.3} We have
\ben
\|\na_h\cdot(\vep\bE_h)\|_\cM\le C\Theta^{1/2}(\|\al^{1/2}\tau^{-1/2}\lj\vep\bE_h\cdot\bn\rj\|_{{\cF}}
+p^{1/2}\|\varphi_h\|_{M_p(\cM)}),
\een
\revt{where $\Theta=\max_{K\in\cM}\Theta_K$.}
\end{lem}

\begin{proof} By \eqref{xx2} and \eqref{g2} we have
\beq\label{g7}
-(\na_h\cdot(\vep\bE_h),q)_\cM+\la\lj\vep\bE_h\cdot\bn\rj,q^+\ra_{{\cF}}+c(\varphi_h,q)=0.
\eeq
By taking $q=\na_h\cdot(\vep\bE_h)$ and using Lemma \ref{lem:3.4} and \eqref{e6}, we obtain
\ben
& &\|\na_h\cdot(\vep\bE_h)\|_\cM^2\\
&\le&C\|\al^{1/2}\tau^{-1/2}\lj\vep\bE_h\cdot\bn\rj\|_\cF\|\na_h\cdot(\vep\bE_h)\|_\cM+\|\varphi_h\|_{M_p(\cM)}\|\na_h\cdot(\vep\bE_h)\|_{M_p(\cM)}\\
&\le&C\Theta^{1/2}(\|\al^{1/2}\tau^{-1/2}\lj\vep\bE_h\cdot\bn\rj\|_\cF+p^{1/2}\|\varphi_h\|_{M_p(\cM)})\|\na_h\cdot(\vep\bE_h)\|_\cM.
\een
This completes the proof.
\end{proof}

The following lemma shows the consistency error for $\div(\vep\bE-\vep\bE_h)$.

\begin{lem}\label{lem:4.4} For any $\psi\in H^1_0(\Om)$, we have
\ben
|(\vep(\bE-\bE_h),\na\psi)|\le C\Theta^{1/2}\,\frac{h}{p^{1/2}}\|\psi\|_{H^1(\Om)}(\|\bE-\bE_h\|_{\bX_p(\cM)}+\|\vp_h\|_{M_p(\cM)}).
\een
\end{lem}

\begin{proof} For any $\psi\in H^1_0(\Om)$, let $\pi_h(\psi)|_K=\pi_{hp}^K(\psi)\ \ \forall K\in\cM$, where $\pi_{hp}^K:H^1(K)\to Q_p(K)$ is defined in Lemma \ref{lem:4.1}.
By \eqref{x3}, we have
\ben
(\vep(\bE-\bE_h),\na\psi)&=&b(\bE-\bE_h,\psi)\\
&=&b(\bE-\bE_h,\psi-\pi_h(\psi))+c(\vp_h,\pi_h(\psi)).
\een
Since $\div(\vep\bE)=0$ in $\Om$, $\lj\vep\bE\cdot\bn\rj|_F=0\ \ \forall F\in\cF$, by \eqref{xx2}, Lemma \ref{lem:4.3}, and Lemma \ref{lem:4.1} we have
\ben
|(\vep(\bE-\bE_h),\na\psi)|&\le&C\frac hp\|\psi\|_{H^1(\Om)}(\|\na_h\cdot(\vep\bE_h)\|_\cM+\|\al^{1/2}\tau^{-1/2}\lj\vep\bE_h\cdot\bn\rj\|_\cF)\\
& &\,+\,\|\vp_h\|_{M_p(\cM)}\|\pi_h(\psi)\|_{M_p(\cM)}.
\een
\rev{For any $K\in\cM$, by Lemma \ref{lem:4.1}, $\|\tau^{1/2}\tau_1^{1/2}\na\pi_h\psi\|_{L^2(K)}\le Ch_Kp^{-3/2}\|\psi\|_{H^1(K)}$. For any $F\in\overline{\cF}$, since $\lj\psi\rj_F=0$ on $F$, by \eqref{g5}, $\|\tau^{1/2}\lj\pi_h\psi\rj\|_{L^2(F)}=\|\tau^{1/2}\lj\pi_h\psi-\psi\rj\|_{L^2(F)}\le Ch_Fp^{-3/2}\|\psi\|_{H^1(\omega(F))}$, where $\omega(F)$ is the  union of the elements having $F$ as one of its faces. This implies $\|\pi_h\psi\|_{M_p(\cM)}\le Chp^{-3/2}\|\psi\|_{H^1(\Om)}$ and thus
\ben
|(\vep(\bE-\bE_h),\na\psi)|\le C\Theta^{1/2}\,\frac{h}{p^{1/2}}\|\psi\|_{H^1(\Om)}(\|\bE-\bE_h\|_{\bX_p(\cM)}+\|\vp_h\|_{M_p(\cM)}).
\een}
This completes the proof.
\end{proof}

\begin{thm}\label{thm:4.2}
Let $\bE\in\bH^{s}(\Om_1\cup\Om_2)$, $s\ge 2$. If $kh/p^{1/2}$ is sufficiently small, then the problem \eqref{g1}-\eqref{g2} has a unique solution $(\bE_h,\vp_h)\in\bX_p(\cM)\times M_p(\cM)$ \rev{which} satisfies the following error estimate
\beq\label{xx}
\|\bE-\bE_h\|_{\bX_p(\cM)}+\|\vp_h\|_{M_p(\cM)}
\le C\revt{\max_{K\in\cM}\left(\frac{\Theta_K^{1/2}h_K^{\nu-1}}{p^{s-3/2}}\right)}\|\bE\|_{H^s(\Om_1\cup\Om_2)}.
\eeq
\end{thm}

\begin{proof} We only need to prove the estimate \eqref{xx} under the condition of the theorem. The uniqueness follows directly from the estimate \eqref{xx} since from $\bJ=0$ in $\Om$ and $\bg=0$ on $\Sigma$, we have $\bE=0$ in $\Om$ and thus $(\bE_h,\vp_h)=(0,0)$ by \eqref{xx}. The existence is a consequence of the uniqueness.

Now we show the estimate \eqref{xx}. Denote by $(\tilde\bE_h,\tilde\vp_h)=\Pi_h(\bE,0)$, where $\Pi_h$ is the projection operator for the coercive Maxwell equations defined in \eqref{g3}-\eqref{g4}. Then from \eqref{g3}-\eqref{g4} and
\eqref{x2}-\eqref{x3} we have
\be
& &a(\tilde\bE_h-\bE_h,\bv)-\overline{b(\bv,\tilde\vp_h-\vp_h)}+k^2(\vep(\tilde\bE_h-\bE_h),\bv)_\cM\nn\\
& &\hskip3cm=R(\bE,\bv)+2k^2(\vep(\bE-\bE_h),\bv)_\cM\ \ \ \ \forall\bv\in\bX_p(\cM),\label{y1}\\
& &b(\tilde\bE_h-\bE_h,q)+c(\tilde\vp_h-\vp_h,q)=0\ \ \ \ \forall q\in M_p(\cM).\label{y2}
\ee
By taking $\bv=\tilde\bE_h-\bE_h$ in \eqref{y1}, $q=\tilde\vp_h-\vp_h$ in \eqref{y2}, we have by Lemma \ref{lem:5.1} that
\be\label{y3}
& &\|\tilde\bE_h-\bE_h\|_{\bX_p(\cM)}+\|\tilde\vp_h-\vp_h\|_{M_p(\cM)}\nn\\
&\le&C\sup_{0\not=\bv\in\bX_p(\cM)}\frac{|R(\bE,\bv)|+|2k^2(\vep(\bE-\bE_h),\bv)_\cM|}{\|\bv\|_{\bX_p(\cM)}}\nn\\
&\le&C\revt{\max_{K\in\cM}\left(\frac{\Theta_K^{1/2}h_K^{\nu-1}}{p^{s-1}}\right)}\|\bE\|_{H^s(\Om_1\cup\Om_2)}+C\,\sup_{0\not=\bv\in\bX_p(\cM)}\frac{|k^2(\vep(\bE-\bE_h),\bv)_\cM|}{\|\bv\|_{\bX_p(\cM)}}.
\ee
For any $\bv\in\bX_p(\cM)$, we define $\psi\in H^1_0(\Om)$ as the solution of the problem
\beq\label{y4}
(\vep\na\psi,\na q)=k(\vep\bv,\na q)\ \ \ \ \forall q\in H^1_0(\Om),
\eeq
Denote $\bw=k\bv-\na\psi$, then
\beq\label{zzz1}
k\bv=\bw+\na\psi,\ \ \div(\vep\bw)=0\ \ \mbox{in }\Om.
\eeq
Obviously, $\|\vep^{1/2}\na\psi\|_{L^2(\Om)}\le k\|\vep^{1/2}\bv\|_\cM$ and $\|\bw\|_{L^2(\Om)}\le C\|k\bv\|_\cM\le C\|\bv\|_{\bX_p(\cM)}$.

Combining \eqref{y3}, Lemma \ref{lem:4.4}, and Theorem \ref{thm:4.1} for estimating $\|\bE-\tilde{\bE}_h\|_{\bX_p(\cM)}+\|\tilde\vp_h\|_{M_p(\cM)}$, we have by the triangle inequality that
\be
& &\|\bE-\bE_h\|_{\bX_p(\cM)}+\|\vp_h\|_{M_p(\cM)}\nn\\
&\le&C\revt{\max_{K\in\cM}\left(\frac{\Theta_K^{1/2}h_K^{\nu-1}}{p^{s-3/2}}\right)}\|\bE\|_{H^s(\Om_1\cup\Om_2)}+C\Theta^{1/2}\,\frac{kh}{p^{1/2}}(\|\bE-\bE_h\|_{\bX_p(\cM)}+\|\vp_h\|_{M_p(\cM)})\nn\\
& &\ +\,C\sup_{0\not=\bv\in\bX_p(\cM)}\frac{|k(\vep(\bE-\bE_h),\bw)_\cM|}{\|\bv\|_{\bX_p(\cM)}}.\label{y6}
\ee

Now we use the duality argument to estimate the last term in \eqref{y6}. Let $\bz$ be the solution of the problem
\be
& &\na\times(\mu^{-1}\na\times\bz)-k^2\vep\bz=\vep\bw,\ \ \div(\vep\bz)=0\ \ \ \ \mbox{in }\Om,\label{y7}\\
& &\lj\bn\times\bz\rj_\Ga=0,\ \ \lj(\mu^{-1}\na\times\bz)\times\bn\rj_\Ga=0,\ \ \lj\vep\bz\cdot\bn\rj_\Ga=0\ \ \ \ \mbox{on }\Ga,\label{y8}\\
& &\bz\times\bn=0\ \ \ \ \mbox{on }\Sigma.\label{y9}
\ee
Since $\div(\vep\bw)=0$ in $\Om$ by \eqref{zzz1}, we know by the regularity Theorem \ref{thm:2.1} that 
\beq\label{k2}
\|\bz\|_{H^2(\Om_1\cup\Om_2)}=\|\bar\bz\|_{H^2(\Om_1\cup\Om_2)}\le C_{\rm reg}\|\vep\bw\|_{L^2(\Om)}\le CC_{\rm reg}\|\bv\|_{\bX_p(\cM)}.
\eeq
Multiplying \eqref{y7} by $\overline{\bE-\bE_h}$ and integrating by parts,
\ben
(\vep(\bE-\bE_h),\bw)_\cM&=&(\mu^{-1}\na_h\times(\bE-\bE_h),\na\times\bz)_\cM-\la\lj\bn\times(\bE-\bE_h)
\rj,\mu^{-1}\na\times\bz\ra_{\overline{\cF}}\\
& &\ -\,k^2(\vep(\bE-\bE_h),\bz)_\cM.
\een
Let $\bz_I=I_{h,p}(\bz)\in\bX_p(\cM)$ be defined in Lemma \ref{lem:4.5}. By \eqref{x2}
\ben
a(\bE-\bE_h,\bz_I)+\overline{b(\bz_I,\vp_h)}-k^2(\vep(\bE-\bE_h),\bz_I)=R(\bE,\bz_I).
\een
Thus
\be\label{k3}
& &(\vep(\bE-\bE_h),\bw)_\cM\nn\\
&=&(\mu^{-1}\na_h\times(\bE-\bE_h),\na_h\times(\bz-\bz_I))_\cM\nn\\
& &\ -\,\la\lj\bn\times(\bE-\bE_h)\rj,\mu^{-1}\na\times\bz-(\mu^{-1}\na_h\times\bz_I)^-\ra_{\overline{\cF}}\nn\\
& &\ -\,k^2(\vep(\bE-\bE_h),\bz-\bz_I)_\cM+(\mu^{-1}\na_h\times(\bE-\bE_h)-\mathsf{L}(\bE-\bE_h),\mathsf{L}(\bz_I))_\cM\nn\\
& &\ -\,\la\al\tau^{-1}\lj\bn\times(\bE-\bE_h)\rj,\lj\bn\times\bz_I\rj\ra_{\overline{\cF}}
-\la\al\tau^{-1}\lj\vep(\bE-\bE_h)\cdot\bn\rj,\lj\vep\bz_I\cdot\bn\rj\ra_{{\cF}}\nn\\
& &\ +\,R(\bE,\bz_I)-\overline{b(\bz_I,\vp_h)}
:={\rm I}_1+\cdots+{\rm I}_8.
\ee
By Lemma \ref{lem:4.5}, we have
\ben
|{\rm I}_1+{\rm I}_2+{\rm I}_3|&\le&C\frac{h}{p}\left(1+\frac{kh}p\right)\|\bz\|_{H^2(\Om_1\cup\Om_2)}\|\bE-\bE_h\|_{\bX_p(\cM)}.
\een
By \eqref{m3} and Lemma \ref{lem:4.5}, 
\ben
|{\rm I}_4+{\rm I}_5+{\rm I}_6|
&\le&C\|\al^{1/2}\tau^{-1/2}\lj\bn\times(\bz-\bz_I)\rj\|_{\overline{\cF}}\|\bE-\bE_h\|_{\bX_p(\cM)}\\
& &\ +\,\|\al^{1/2}\tau^{-1/2}\lj\vep(\bz-\bz_I)\cdot\bn\rj\|_{{\cF}}\|\bE-\bE_h\|_{\bX_p(\cM)}\\
&\le&C\Theta^{1/2}\frac h{p^{1/2}}\|\bz\|_{H^2(\Om_1\cup\Om_2)}\|\bE-\bE_h\|_{\bX_p(\cM)}.
\een
By Lemma \ref{lem:5.1} and Lemma \ref{lem:4.5}, 
\ben
|{\rm I}_7|
&\le&C\revt{\max_{K\in\cM}\left(\frac{\Theta_K^{1/2}h_K^{\nu-1}}{p^{s-1}}\right)}\|\bE\|_{H^s(\Om_1\cup\Om_2)}\|\tau^{-1/2}\lj\bn\times\bz_I\rj\|_\cF\nn\\
&\le&C\revt{\max_{K\in\cM}\left(\frac{\Theta_K^{1/2}h_K^{\nu-1}}{p^{s-1}}\right)}\frac{h}{p^{1/2}}\|z\|_{H^2(\Om_1\cup\Om_2)}\|\bE\|_{H^s(\Om_1\cup\Om_2)}.
\een
Finally, since $\div(\vep\bz)=0$ in $\Om$ and $\lj\vep\bz\cdot\bn\rj_{\Ga}=0$ on $\Ga$, we obtain by Lemma \ref{lem:4.5} that
\ben
|{\rm I}_8|
&\le&(\|\tau^{-1/2}\tau_1^{-1/2}(\bz_I-\bz)\|_\cM+\|\tau^{-1/2}(\bz_I-\bz)\|_{\overline{\cF}})\|\vp_h\|_{M_p(\cM)}\nn\\
&\le&C\frac h{p^{1/2}}\|\bz\|_{H^2(\Om_1\cup\Om_2)}\|\vp_h\|_{M_p(\cM)}.
\een
Inserting the estimates for ${\rm I}_1,\cdots,{\rm I}_8$ into \eqref{k3}, we obtain by \eqref{y6} and \eqref{k2} that
\ben
\|\bE-\bE_h\|_{\bX_p(\cM)}+\|\vp_h\|_{M_p(\cM)}
&\le&C\revt{\max_{K\in\cM}\left(\frac{\Theta_K^{1/2}h_K^{\nu-1}}{p^{s-3/2}}\right)}\|\bE\|_{H^s(\Om_1\cup\Om_2)}\\
& &\,+\,CC_{\rm reg}\Theta^{1/2}\,\frac{kh}{p^{1/2}}(\|\bE-\bE_h\|_{\bX_p(\cM)}+\|\vp_h\|_{M_p(\cM)}).
\een
Therefore, if ${kh}/{p^{1/2}}$ is sufficiently small, the desired estimate \eqref{xx}
in the theorem is proved. This completes the proof.
\end{proof}

%%%%%%%%%%%%%%%%%%%%%%%%%%%%%%%%%%%%%%%%%%%%%%%%%%%
\section{Numerical results}
%%%%%%%%%%%%%%%%%%%%%%%%%%%%%%%%%%%%%%%%%%%%%%%%%%%

In this section, we show a numerical example to illustrate the performance of the proposed unfitted finite element method.
Let $\Om = (-1,1)^3$ and $\Om_1=\{(x_1,x_2,x_3)\in\R^3: \sum^3_{i=1}x_i^2/d_i^2<1\}$, where $d_1=d_2=0.4$ and $d_3=0.8$, and $\Om_2=\Om \backslash \bar{\Om}_1$. Set $\mu_1=\ep_1=2$, $\mu_2=\ep_2=1$, and $k=1$. Functions $\bJ$ and $\bg$ are chosen such that the exact solution $\bE$ to the problem \eqref{p1}-\eqref{p4} is given by 
\ben
\bE(\bx)=\na\times\left[\mu_i\Big(\frac{x_1^2}{d_1^2}+\frac{x_2^2}{d_2^2}+\frac{x_3^2}{d_3^2}-1\Big)^2\cos(10x_1)
\left(\begin{array}{c}
1 \\
1\\
1\\
\end{array}\right)\right],\ \  \bx\in\Om_i,\ \ i=1,2.
\een
%The computations are carried out on a cluster with Intel Xeon Gold 6230.

We start with a uniform Cartesian mesh $\cT$ of $\Om$ with the mesh size $h_0=2/n$ with $n=8,16,32,\cdots$. The induced mesh $\cM$ satisfying Assumption (H1) and the upper bound of the interface deviation $\eta$ is constructed by using a merging algorithm developed by Linbo Zhang in the software platform Parallel Hierarchical Grid (PHG) \cite{PHG}, which iteratively refines the interface elements and their neighboring elements in $\cT$ such that $\eta_K\le\eta\ \forall K\in\cM^\Ga$. The finite element bases for the elements $K\in\cM\backslash\cM^\Ga$ are the Lagrangian interpolation polynomials through the Gauss-Lobatto integration points. The bases for the interface elements $K\in\cM^\Ga$ are $L^2$ orthogonal functions on the maximal polyhedron inside $K_i=K\cap\Om_i$, $i=1,2$, which extend a similar two-dimensional construction in \cite{Ch22}. High-order numerical integration in $K_i$, $i=1,2$, and $\Gamma_K$ for $K\in\mathcal{M}^{\Gamma}$ is performed by using numerical quadrature functions in Cui et al. \cite{Cui}. The resultant linear systems of equations are solved by the preconditioned GMRES method with the overlapping additive Schwarz method as the preconditioner implemented in PHG. We set $\de_0=1/4$ in the definition of the large element in Definition \ref{def:3.1}.

To choose the upper bound of the interface deviation $\eta$, we note that by Theorem \ref{thm:4.2}, the error is of the order $\max_{K\in\cM}(\Theta_K^{1/2}h_K^p)$ for smooth solutions. We know that for the elements $K$ away from the interface, $\Theta_{K}=1$ and $h_{K}=h_0$, and by Lemma \ref{lem:3.0}, $\eta_K\le Ch_K$ for any $K\in\cM^\Ga$. This motivates us to choose $\eta$ by the equation
\ben
h_0^p=\beta\Theta^{1/2}\eta^p\ \ \iff\ \ h_0^p=\beta\mathsf{T}\left(\frac{1+3\eta}{1-\eta}\right)^{3p}\eta^p,
\een 
where $\beta>0$ is some fixed constant. Recall that $\mathsf{T}(t)=t+\sqrt{t^2-1}$ for $t\ge 1$. In our computations we choose $\beta=0.1$.

\begin{figure}[h]
	\begin{center}
		\begin{minipage}{5cm}
			\epsfxsize=0.9\textwidth\epsffile{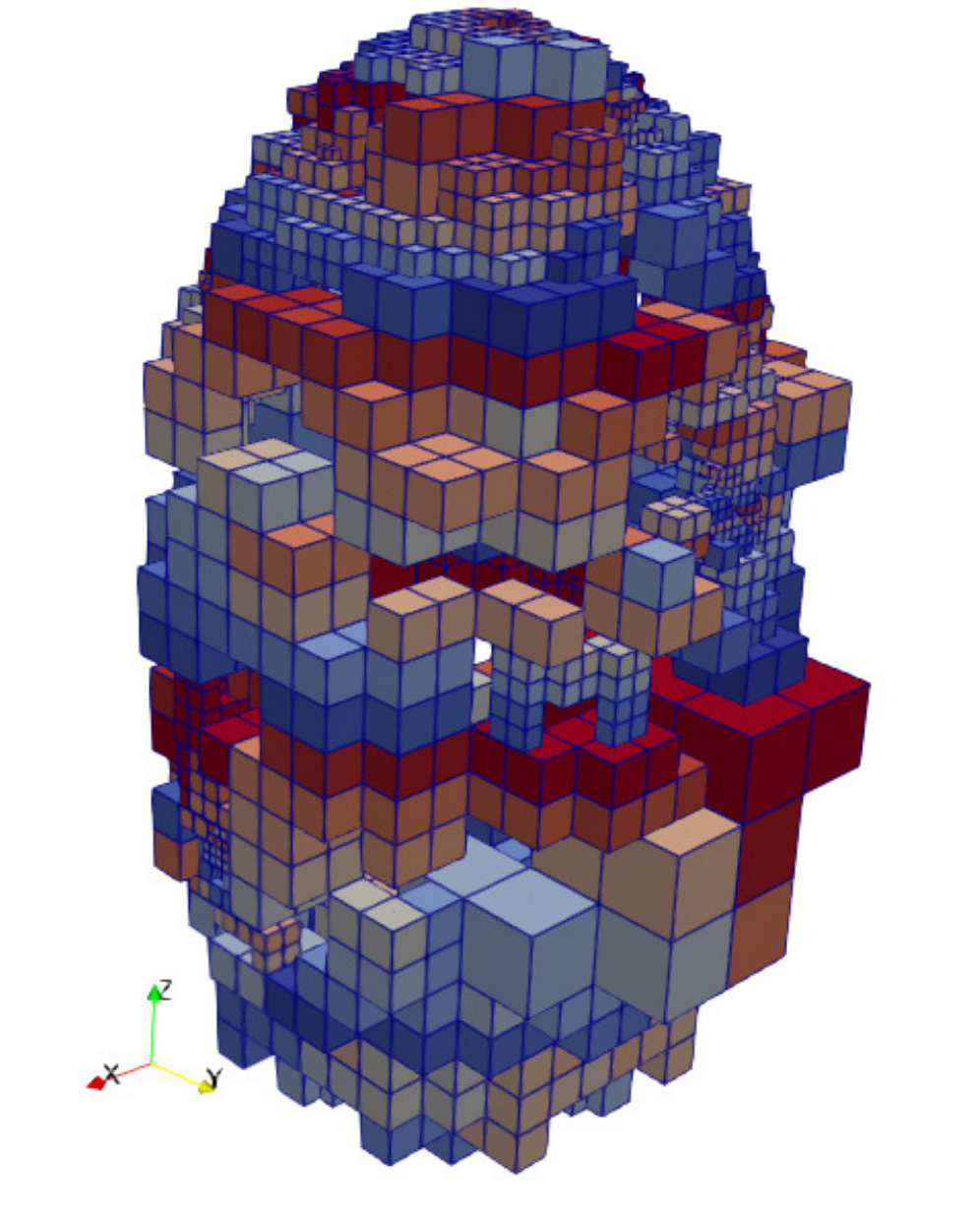}
		\end{minipage}
		\begin{minipage}{3.5cm}
			\epsfxsize=1\textwidth\epsffile{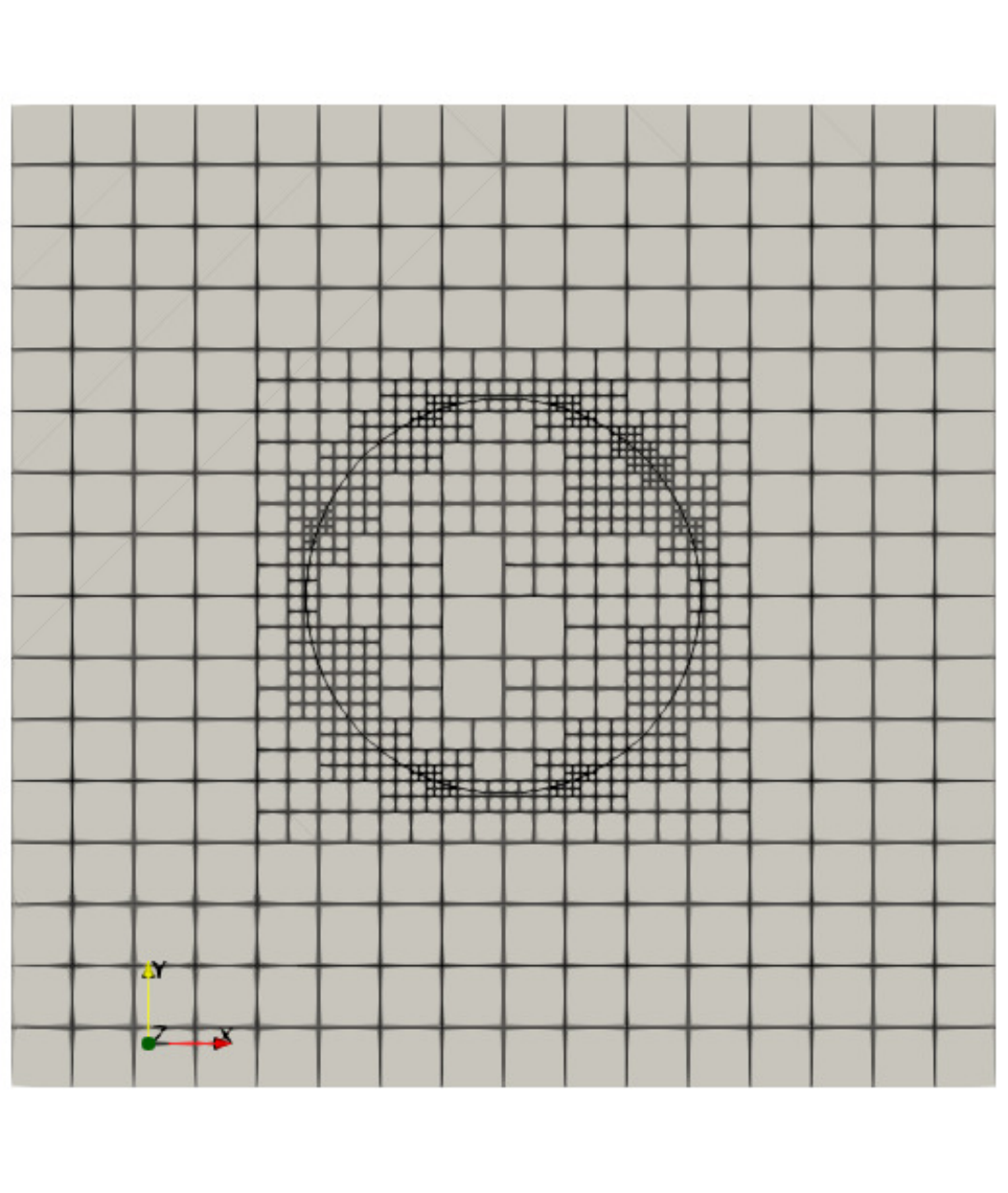}
			\end{minipage}
			\begin{minipage}{3.5cm} %\vskip-0.1cm
			\epsfxsize=1\textwidth\epsffile{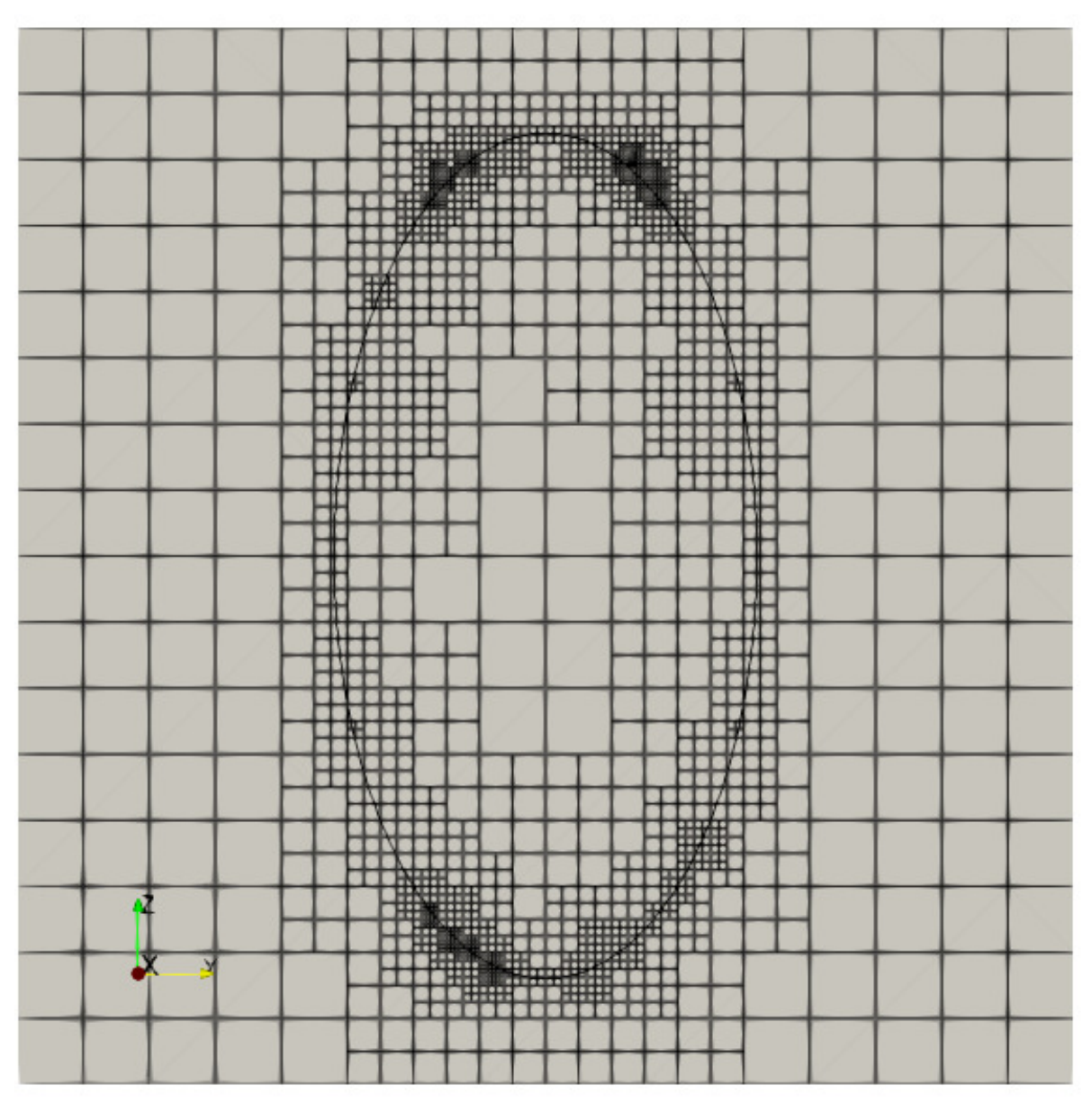}
		\end{minipage}
	\end{center}
	\caption{The induced mesh across $\Gamma$ in $\Omega_1$ when $h_0=1/8, p=1$ (left), the cross-sections of the induced mesh at $x_3=0$ (middle) and $x_1=0$ (right) when $h_0=1/8$, $p=3$.}\label{fig:5.2}
\end{figure}

Fig.\ref{fig:5.2} shows the induced mesh across $\Gamma$ in $\Omega_1$ in which only macro-elements merged by more than 2 elements in $\cT$ are shown, the cross-sections of the induced mesh at $x_3=0$ and $x_1=0$ when $h_0=1/8$. Table \ref{tab:5.1} shows the upper bound of the interface deviation $\eta$, the number of elements of the induced mesh $N$, the number of the interface elements $N^\Ga$, the number of degrees of freedom $\mbox{\#DoFs}(\bX_p)$ and $\mbox{\#DoFs}(M_p)$ for $\bX_p(\cM)$ and $M_p(\cM)$, the relative error $\mathcal{E}=(\|\bE-\bE_h\|_{\bX_p(\cM)}+\|\vp_h\|_{M_p(\cM)})/\|\bE\|_{\bX_p(\cM)}$, and the convergence order of the $p$-th order methods when $p=1,2,3$. We observe the optimal convergence of the relative error $\mathcal{E}$ which is in conform with Theorem \ref{thm:4.2}. We also observe that to achieve the same relative error, high order method is much more efficient by requiring less number of degrees of freedom. This shows clearly the advantage of using high order methods. We remark that \#DoFs$(\bX_p)$ and \#DoFs$(M_p)$ can be further reduced if one uses conforming finite elements away from the interface, which will be implemented in our code in future.

\begin{table}[htb]
\caption{\rm The optimal convergence of the unfitted finite element method when $p=1,2,3$.}
  \centering
  \fontsize{9pt}{12pt}\selectfont
    \vspace{0.1in}
    \begin{tabular}{cccccccc}
    \hline
	$h$  
	& $\eta$  
    & $N$
    & $N^\Gamma$
    & {\#DoFs($\bX_p$)} & {\#DoFs($M_p$)}
    & $\mathcal{E}$ & order
    \cr\cline{1-8}
    \multicolumn{8}{c}{$p=1$}
    \cr\cline{1-8}    
    $1/4$ & 1.21e-01 &    5,412 &  1,162 &      157,776 &     52,592 & 5.81e-01 &     --
    \cr\cline{1-8}
    $1/8$ & 9.17e-02 &   14,631 &  2,773 &     417,696 &    139,232 & 2.27e-01 &   1.36
    \cr\cline{1-8}
   $1/16$ & 6.73e-02 &   46,978 &  4,981 &    1,247,016 &    415,672 & 1.14e-01 &   0.99
    \cr\cline{1-8}
   $1/32$ & 4.80e-02 &  291,789 & 12,524 &   12,722,712 &  4,240,904 & 5.66e-02 &   1.01
    \cr\cline{1-8}   
   $1/64$ & 3.31e-02 &2,151,108 & 34,363 &   52,451,304 & 17,483,768 & 2.84e-02 &   0.99
    \cr\cline{1-8}    
    \multicolumn{8}{c}{$p=2$}
    \cr\cline{1-8}  
    $1/4$ & 7.50e-02 &   20,238 &  4,331 &    19,90,089 &   6,63,363 & 1.47e-01 &     -- 
    \cr\cline{1-8}
    $1/8$ & 5.40e-02 &   49,715 & 10,625 &    4,887,540 &  1,629,180 & 3.83e-02 &   1.94
    \cr\cline{1-8}
   $1/16$ & 3.77e-02 &  111,973 & 19,842 &   10,677,015 &  3,559,005 & 9.65e-03 &   1.99
    \cr\cline{1-8}
   $1/32$ & 2.54e-02 &  440,378 & 47,022 &   39,479,400 & 13,159,800 & 2.51e-03 &   1.94
    \cr\cline{1-8}
   $1/64$ & 1.65e-02 &2,508,941 &116,157 &  212,632,938 & 70,877,646 & 6.48e-04 &   1.95
    \cr\cline{1-8}  
    \multicolumn{8}{c}{$p=3$}
    \cr\cline{1-8}  
    $1/4$ & 6.28e-02 &   30,073 &  6,357 &    6,994,560 &  2,331,520 & 3.84e-02 &     --  
    \cr\cline{1-8}
    $1/8$ & 4.44e-02 &   68,167 & 15,132 &   15,993,408 &  5,331,136 & 5.02e-03 &   2.94
    \cr\cline{1-8}
   $1/16$ & 3.04e-02 &  172,040 & 32,484 &   39,268,608 & 13,089,536 & 7.25e-04 &   2.79
    \cr\cline{1-8}
   $1/32$ & 2.01e-02 &  546,078 & 71,425 &  118,560,576 & 39,520,192 & 9.12e-05 &   2.99
    \cr\cline{1-8}   
\end{tabular}\label{tab:5.1}
\end{table}

\medskip

%%%%%%%%%%%%%%%%%%%%%%%%%%%%%%%%%%%%%%%%%%%%%%%%%%%
{\bf Acknowledgement.} The authors are very thankful to the referees for the constructive comments that lead to great improvement of the paper. The authors would also like to
thank Prof. Linbo Zhang and Dr. Yong Liu from Chinese Academy of Sciences for the inspiring discussion and generous
help in performing the numerical experiments. 
%%%%%%%%%%%%%%%%%%%%%%%%%%%%%%%%%%%%%%%%%%%%%%%%%%%

\end{document}